\documentclass[12pt]{article}
\usepackage{amsmath, amsfonts, amssymb}
\newcommand{\dal}{\square}
\newcommand{\R}{{\mathbb R}}
\newcommand{\C}{{\mathbb C}}
\newcommand{\ve}{\varepsilon}
\newcommand{\pa}{\partial}
\newcommand{\jb}[1]{\left\langle #1 \right\rangle}
\newcommand{\Sum}{\sideset{}{'}\sum}
\DeclareMathOperator{\supp}{{\rm supp}}
\DeclareMathOperator{\diag}{{\rm diag}}
\DeclareMathOperator{\RealPart}{{\rm Re}}
\DeclareMathOperator{\ImPart}{{\rm Im}}
\newtheorem{theorem}{Theorem}
\newtheorem{lemma}{Lemma}
\newtheorem{proposition}{Proposition}
\newtheorem{definition}{Definition}
\newtheorem{example}{Example}
\numberwithin{equation}{section}
\numberwithin{theorem}{section}
\numberwithin{lemma}{section}
\numberwithin{proposition}{section}
\numberwithin{definition}{section}
\begin{document}
\title{%
Global Existence for Coupled Systems of Nonlinear Wave and Klein-Gordon Equations
in Three Space Dimensions
}
\author{Soichiro Katayama}
\date{}
%\institute{S. Katayama \at Department of Mathematics, Wakayama University, 930 Sakaedani, Wakayama 640-8510, Japan
%\\ Tel. +81-73-457-7343, Fax: +81-73-457-7489
%\\ \email{katayama@center.wakayama-u.ac.jp}
%}
%\titlerunning{Systems of Nonlinear Wave and Klein-Gordon Equations}
%%%%%%
\maketitle
%%%%%%%%%%%%
%%%%%%%%%%%%
\begin{abstract}
We consider the Cauchy problem for coupled systems of wave and Klein-Gordon equations
with quadratic nonlinearity in three space dimensions. 
We show global existence of small amplitude solutions
under certain condition including 
the null condition on self-interactions 
between wave equations.
Our condition is much weaker than
the strong null condition introduced by Georgiev for this kind of coupled system.
Consequently our result is applicable to certain physical systems,
such as the Dirac-Klein-Gordon equations, the Dirac-Proca equations,
and the Klein-Gordon-Zakharov equations. 
\end{abstract}
%%%%%%%%%%%%
%\keywords{Null condition \and Wave equation \and Klein-Gordon equation \and Global existence}
%\subclass{35L70}
%%%%%%%%%%%%
\section{Introduction}
We consider the Cauchy problem for the following system:
\begin{align}
\label{OurSystem}
\left(\dal+m_i^2\right) u_i=F_i(u, \pa u, \pa_x\pa u), \quad i=1,2, \ldots, N
\end{align}
in $(0, \infty)\times \R^3$ with initial data
\begin{equation}
\label{Data} u(0,x)=\ve f(x),\ (\pa_t u)(0, x)=\ve g(x) \text{ for $x=(x_1, x_2, x_3)\in \R^3$},
\end{equation}
where $\dal=\pa_t^2-\Delta_x$, $u=(u_j)_{1\le j\le N}$, $m_i\ge 0$ for $1\le i\le N$,
and $\varepsilon$ is a small and positive parameter. Here 
each component $u_j$ of $u$ is supposed to be a real-valued unknown function of $(t,x)\in (0,\infty)\times \R^3$,
and $\Delta_x$ denotes the Laplacian in $x$-variables.
In the above, $\pa u$ and $\pa_x\pa u$ are given by
$$
\pa u=(\pa_a u_j)_{1\le j\le N, 0\le a\le 3}, \quad 
\pa_x\pa u=(\pa_k\pa_a u_j)_{1\le j\le N, 1\le k\le 3, 0\le a\le 3},
$$
respectively, with the notation 
$$
\pa_0=\pa_t=\frac{\pa}{\pa t} \text{ and } \pa_k=\frac{\pa}{\pa x_k} \text{ for }1\le k\le 3.
$$
Here, by writing $(\pa_a u_j)_{j, a}$, 
we mean that $\pa_a u_j$'s are arranged in dictionary order with respect to $(j,a)$.
Similarly, $(\pa_k \pa_ a u_j)_{j,k,a}$ means that 
$\pa_k\pa_a u_j$'s are arranged in dictionary order with respect to $(j,k,a)$.
Similar convention will be used throughout this paper.
We assume that, for $1\le i\le N$, each $F_i=F_i(\xi, \xi', \xi'')$ is
a real-valued smooth function of $(\xi, \xi', \xi'')\in \R^N\times \R^{4N}\times \R^{12N}$,
where $\xi$, $\xi'$, and $\xi''$ are independent variables for which
$u$, $\pa u$, and $\pa_x\pa u$ are substituted in \eqref{OurSystem};
more precisely, if we write 
\begin{align*}
& \xi=(\xi_j)_{1\le j\le N}\in \R^N, \quad \xi'=(\xi_{j, a}')_{1\le j\le N, 0\le a\le 3}\in \R^{4N},\\
& \xi''=(\xi_{j, k, a}'')_{1\le j\le N, 1\le k\le 3, 0\le a\le 3} \in \R^{12N},
\end{align*}
then $\xi_j$, $\xi_{j,a}'$, and $\xi_{j,k,a}''$ are the independent variables 
for which 
$u_j$, $\pa_a u_j$, and $\pa_k \pa_a u_j$
are substituted in \eqref{OurSystem}, respectively. 
We assume that $F=(F_i)_{1\le i\le N}$ vanishes of second order at the origin $(\xi, \xi', \xi'')=(0,0,0)$, namely
$$
F(\xi, \xi', \xi'')=O(|\xi|^2+|\xi'|^2+|\xi''|^2) \text{ around $(\xi, \xi', \xi'')=(0,0,0)$}.
$$
%$F=(F_j)_{1\le j\le N}$ is a smooth and $\R^N$-valued function vanishing of second order at the origin.
For simplicity, we also assume that the system is quasi-linear. In other words, we assume that
\begin{equation}
\label{quasi-linear}
F_i(\xi, \xi', \xi'')=\sum_{j=1}^N \sum_{k=1}^3 \sum_{a=0}^3
\gamma_{ka}^{\, ij}(\xi, \xi') \xi''_{j,k,a}+{G}_i(\xi, \xi'), \quad 1\le i\le N,
\end{equation}
%for $1\le i\le N$, 
where $\gamma_{ka}^{\, ij}=\gamma_{ka}^{\, ij}(\xi, \xi')$
and ${G}_i={G}_i(\xi, \xi')$ are some functions vanishing of first and second order
at the origin, respectively.
Moreover, to assure the hyperbolicity of the system, 
we always assume the symmetricity condition
\begin{equation}
\label{symmetricity}
\gamma_{ka}^{\, ij}(\xi, \xi')=\gamma_{ka}^{\, ji}(\xi, \xi') \text{ and }
\gamma_{kl}^{\, ij}(\xi, \xi')=\gamma_{l k}^{\, ij}(\xi, \xi')
\end{equation}
for any $(\xi, \xi')\in \R^N\times \R^{4N}$, $1\le i,j \le N$, $1\le k, l \le 3$, and $0\le a\le 3$ (note that the last half
of \eqref{symmetricity} is no restriction as far as we consider smooth solutions).
%%%%%%%%%%%%%%%%%%
For each $i\in \{1, \ldots, N\}$, the equation \eqref{OurSystem} for $u_i$ is called
a (nonlinear) {Klein-Gordon} (resp.~wave) {equation} if $m_i>0$ (resp.~$m_i=0$).

For a while, we suppose that $f=(f_i)_{1\le i\le N}$, $g=(g_i)_{1\le i\le N}\in C^\infty_0(\R^3; \R^N)$ in \eqref{Data}.
Under the conditions \eqref{quasi-linear} and \eqref{symmetricity}, 
the classical theory of nonlinear hyperbolic systems implies
local existence of smooth solutions to the Cauchy problem \eqref{OurSystem}--\eqref{Data}
for sufficiently small $\ve$. Hence we are interested in the sufficient condition for global existence of small amplitude solutions. 
Here we recall the known results briefly.
If the nonlinearity $F$ vanishes of third order at the origin,
\eqref{OurSystem}--\eqref{Data} admits a global solution for small $\ve$.
For arbitrary quadratic nonlinearity $F$, we also have global existence
of small solutions if \eqref{OurSystem} is a system
of nonlinear Klein-Gordon equations, 
namely if $m_i>0$ for all $i=1, \ldots, N$
(see Klainerman~\cite{Kla85:02} and Shatah~\cite{Sha85}; see also Bachelot~\cite{Bac88} and
Hayashi-Naumkin-Ratno Bagus Edy Wibowo~\cite{Hay-Nau-Bag08}).
By contrast, this is not true if
\eqref{OurSystem} is
a system of wave equations (namely
if $m_1=m_2=\cdots =m_N=0$), and the solution to \eqref{OurSystem}--\eqref{Data} with certain
quadratic nonlinearity $F$ may blow up in finite time
no matter how small $\ve$ is (see John \cite{Joh79}, \cite{Joh81}).
Thus we need to put some condition on quadratic nonlinearity,
in order to obtain global solutions for wave equations.
The {\it null condition} introduced by 
Klainerman~\cite{Kla86} is one of such conditions.
%%%%%%%%%%%%%%%%%%%%%%%%%
Before describing the null condition, we introduce the following notation:
For a smooth function  $\Phi=\Phi(z)$ ($z\in \R^d$), we write $\Phi^{\rm (q)}$ for
the quadratic part of $\Phi$ (``(q)'' stands for ``quadratic''); more precisely,
for a smooth function $\Phi=\Phi(z)$, we define
\begin{equation}
\label{qp}
\Phi^{\rm (q)}(z)=\sum_{|\alpha|=2} \frac{(\pa_z^\alpha \Phi)(0)}{\alpha!}z^\alpha,\quad z=(z_1, \ldots, z_d)\in \R^d,
\end{equation}
where $\pa_z=(\pa_{z_1}, \ldots, \pa_{z_d})$, $\alpha$ is a multi-index, and we have used the 
standard notation of multi-indices.
The null condition can be stated as follows:
\begin{definition}[The null condition]
\label{NullD}
We say that a function $F=(F_i)_{1\le i\le N}$ of 
%%%%%%%%%
%%%%%%%%%
$(\xi, \xi', \xi'')\in \R^N\times\R^{4N}\times\R^{12N}$
satisfies the null condition if each $F_i$ $(1\le i \le N)$ 
satisfies
\begin{equation}
\label{NullC}
F_i^{\rm (q)}\bigl(\lambda, (X_a \mu_j)_{1\le j\le N, 0\le a\le 3}, (X_k X_a \nu_j)_{1\le j\le N, 1\le k\le 3, 0\le a\le 3}\bigr)=0
\end{equation}
for any $\lambda=(\lambda_j)_{1\le j\le N}$, $\mu=(\mu_j)_{1\le j\le N}$, $\nu=(\nu_j)_{1\le j\le N}\in \R^N$,
and any
$X=(X_a)_{0\le a\le 3}$ $\in \R^4$ satisfying $X_0^2-X_1^2-X_2^2-X_3^2=0$,
where $F_i^{\rm (q)}=F_i^{\rm (q)}(\xi, \xi', \xi'')$ is the quadratic part of $F_i$ given by \eqref{qp} with 
$\Phi=F_i$ and $z=(\xi, \xi', \xi'')$.
\end{definition}
If $F$ satisfies the null condition, then we have global existence of small solutions 
for systems of wave equations \eqref{OurSystem} with $m_1=\cdots=m_N=0$
(see Klainerman~\cite{Kla86} and Christodoulou~\cite{Chr86}).

Klainerman used the so-called vector field method 
in \cite{Kla85:02} and \cite{Kla86}.
But his method is not directly applicable to systems 
consisting of both wave and Klein-Gordon equations,
because the scaling operator
$S=t\pa_t+\sum_{k=1}^3x_k \pa_k$,
which was used in \cite{Kla86},
is compatible with 
the wave equations, but not with the Klein-Gordon equations. 
This causes some difficulty in the treatment of the null condition,
and hence Georgiev~\cite{Geo90} introduced the {\it strong null condition} to obtain global existence
of small solutions for coupled systems of nonlinear wave and Klein-Gordon equations
(see Section 4 below for the detail),
where $F$ is said to satisfy the strong null condition if \eqref{NullC} with $1\le i\le N$ holds
for any $\lambda$, $\mu$, $\nu \in \R^N$ and any $X\in \R^4$
{\it not} necessarily satisfying $X_0^2-X_1^2-X_2^2-X_3^2=0$.

Our aim in this paper is to establish
a global existence theorem for systems
of the nonlinear wave and Klein-Gordon equations under more
natural and weaker condition than the strong null condition,
so that it can cover the previous results for
wave equations and the Klein-Gordon equations,
as well as some important examples from physics.
%%%%%%%%%%%%%%%%%%%%%%%%%%%%%%%%%%%
\section{The Main Result and Examples}
%%%%%%%%%%%%%%%%%%%%%%%%%%%%%%%%%%%
First we introduce some notation.
Suppose that we can take some natural number $N_1$ such that we have
\begin{equation}
\label{MasslessMass}
m_i>0\text{ for $1\le i\le N_1$, and } m_i=0 \text{ for $N_1+1\le i\le N$}
\end{equation}
in \eqref{OurSystem}.
We set 
\begin{equation}
\label{division}
v=(v_j)_{1\le j\le N_1}:= (u_j)_{1\le j\le N_1} \text{ and }
w=(w_j)_{1\le j\le N_2} := (u_{N_1+j})_{1\le j\le N_2},
\end{equation}
where $N_2=N-N_1$,
so that $u=(u_1, \ldots, u_{N_1}, u_{N_1+1}, \ldots, u_N)=(v, w)$.
Note that each $v_j(=u_j)$ satisfies a nonlinear Klein-Gordon equation,  
while each $w_j(=u_{N_1+j})$ satisfies a nonlinear wave equation. 
In accordance with \eqref{division},
we introduce independent variables $(\eta, \zeta)\in \R^{N_1}\times \R^{N_2}$,
 $(\eta', \zeta')\in \R^{4N_1}\times \R^{4N_2}$, and $(\eta'', \zeta'')\in \R^{12N_1}\times \R^{12N_2}$ to write
\begin{align*}
\xi&= (\xi_j)_{1\le j\le N}=:\left((\eta_j)_{1\le j\le N_1}, (\zeta_j)_{1\le j\le N_2}\right)=(\eta, \zeta), \\
\xi'&= (\xi'_{j,a})_{1\le j\le N, 0\le a\le 3}
=:\left((\eta'_{j,a})_{1\le j\le N_1,0\le a\le 3}, (\zeta'_{j,a})_{1\le j\le N_2, 0\le a\le 3}\right)=(\eta', \zeta'), \\
\xi''&= (\xi''_{j, k, a})_{1\le j\le N, 1\le k\le 3, 0\le a\le 3}\\
    &=:\left((\eta''_{j, k, a})_{1\le j\le N_1, 1\le k\le 3, 0\le a\le 3}, (\zeta''_{j, k, a})_{1\le j\le N_2, 1\le k\le 3, 0\le a\le 3}\right)=(\eta'', \zeta'').
\end{align*}
Correspondingly, we write $\pa u=(\pa v, \pa w)$ and 
$\pa_x\pa u=(\pa_x\pa v, \pa_x\pa w)$.
For a smooth function $\Phi=\Phi(\xi,\xi',\xi'')$, we define
\begin{align}
\label{FIW}
\Phi^{\rm (W)}(\zeta, \zeta', \zeta'')=& 
        \Phi^{\rm (q)}\bigl((\eta, \zeta), (\eta', \zeta'), (\eta'', \zeta'') \bigr)\Bigr|_{(\eta, \eta', \eta'')=(0,0,0)}
\end{align}
for $(\zeta, \zeta',\zeta'')\in \R^{N_2}\times \R^{4N_2}\times \R^{12N_2}$, 
where $\Phi^{\rm (q)}$ is the quadratic part of $\Phi$ given by \eqref{qp} with $z=(\xi, \xi', \xi'')$. Thus $F_i^{\rm (W)}(w, \pa w, \pa_x\pa w)$ appearing below
represents the self-interaction between the solutions to wave equations
(``(W)'' in the notation \eqref{FIW} stands for ``wave''). 

Now we are in a position to state our main result.
%%%%%%%%%%%%%%%%%%%%%%%%%%%%%%%%
\begin{theorem}\label{Global}
Suppose that \eqref{quasi-linear}, \eqref{symmetricity}, and
\eqref{MasslessMass} are fulfilled.
Assume that the following two conditions {\rm (a)} and {\rm (b)} hold:
\begin{enumerate}
%%%%%%%%%%%%%%%%%%%%%%%%%%%%%%%%
\item[{\rm (a)}] 
$\bigl(F_i^{\rm (W)} \bigr)_{N_1+1\le i\le N}$ satisfies the null condition.
%%%%%%%%%%%%%%%%%%%%%%%%%%%%%%%%
\item[{\rm (b)}] 
There exist two (empty or non-empty) sets ${\mathcal I}_1$ and ${\mathcal I}_2$
satisfying 
\begin{equation}
{\mathcal I}_1\cup {\mathcal I_2}=\{1, \ldots, N_2\}, \quad {\mathcal I}_1\cap {\mathcal I}_2=\emptyset,
\end{equation}
and the following properties:
\begin{itemize}
\item[{\rm (b--i)}] 
For any $k\in {\mathcal I}_1$, we have 
\begin{equation}
\label{b-i-1}
\frac{\pa F_i^{\rm (q)}}{\pa \zeta_k}\bigl((\eta, \zeta), (\eta',\zeta'), (\eta'',\zeta'')\bigr)
\left(=\frac{\pa F_i^{\rm (q)}}{\pa \xi_{N_1+k}}(\xi,\xi',\xi'')\right)=0
\end{equation}
for all $i=1, \ldots, N$, and all $(\xi,\xi',\xi'')\in \R^N\times \R^{4N}\times \R^{12N}$.
%%%%%%%%%%%%%%%%%% 
\item[{\rm (b--ii)}] 
For any $k\in {\mathcal I}_2$, there exist some functions 
${\mathcal G}_{k,a}={\mathcal G}_{k,a}(\xi, \xi')$ with $0\le a\le 3$ such that 
\begin{equation}
F_{N_1+k}\bigl(\phi, \pa \phi, \pa_x \pa \phi \bigr)
=\sum_{a=0}^3 \pa_a \left\{ {\mathcal G}_{k, a}\bigl(\phi, \pa\phi \bigr)\right\}
\label{b-ii-1}
\end{equation}
holds for any $\phi=\phi(t,x)\in C^2\left((0,\infty)\times \R^3; \R^N\right)$,
and 
\begin{equation}
\label{b-ii-2}
\frac{\pa {\mathcal G}_{k, a}^{\rm (q)}}{\pa\zeta_l} \bigl((\eta, \zeta), (\eta', \zeta')\bigr)
\left(=\frac{\pa {\mathcal G}_{k, a}^{\rm (q)}}{\pa\xi_{N_1+l}} (\xi, \xi') \right)=0,\quad 0\le a\le 3
%\equiv 0
\end{equation}
holds for all $l\in {\mathcal I}_1$, and all $(\xi, \xi')\in \R^N\times \R^{4N}$,
where ${\mathcal G}_{k,a}^{\rm (q)}$ is the quadratic part of ${\mathcal G}_{k,a}$.
 \end{itemize}
%%%%%%%%%%%%%%%%%%%%%%
\end{enumerate}
Then, for any $f$, $g\in {\mathcal S}(\R^3; \R^N)$,
there exists a positive constant $\ve_0$ such that
for any $\ve\in (0,\ve_0]$
the Cauchy problem \eqref{OurSystem}--\eqref{Data}
admits a unique global solution $u\in C^\infty\bigl([0, \infty)\times \R^3; \R^N \bigr)$.
Here $\mathcal S$ denotes the Schwartz class, the class of rapidly decreasing functions.
\end{theorem}
%%%%%%%%%%%%%%%%%%%%%%%%%%%%%%%
Here and hereafter, we say that 
$\bigl(F^{\rm (W)}_i\bigr)_{N_1+1\le i\le N}$ satisfies the null condition,
if each $F_i^{\rm (W)}=F_i^{\rm (W)}(\zeta, \zeta', \zeta'')$ with $N_1+1\le i\le N$ 
satisfies
$$
F_i^{\rm (W)} \left(\lambda, (X_a \mu_j)_{1\le j\le N_2, 0\le a\le 3}, (X_kX_a \nu_j)_{1\le j\le N_2, 1\le k\le 3, 0\le a\le 3}\right)=0
$$
for any $\lambda=(\lambda_j)_{1\le j\le N_2}$,
$\mu=(\mu_j)_{1\le j\le N_2}$, $\nu=(\nu_j)_{1\le j\le N_2}\in \R^{N_2}$,
and for any $X=(X_a)_{0\le a\le 3}\in \R^4$ satisfying $X_0^2-X_1^2-X_2^2-X_3^2=0$.
Notice that the null structure is required only for $F_i^{\rm (W)}$ with $N_1+1\le i\le N$ in
Theorem \ref{Global}.
%%%%%%%%%%
We define the {\it null forms} $Q_0$ and $Q_{ab}$ by
\begin{align}
Q_0(\varphi, \psi)= & (\pa_t\varphi)(\pa_t\psi)-(\nabla_x \varphi)\cdot(\nabla_x \psi),\\
Q_{ab}(\varphi, \psi) = & (\pa_a \varphi)(\pa_b \psi)-(\pa_b \varphi)(\pa_a \psi), \quad
0\le a<b\le 3,
\end{align}
where $\nabla_x=(\pa_1, \pa_2, \pa_3)$, and $\cdot$ denotes the inner product in $\R^3$.
Then we can easily check 
that the assumption (a) in Theorem \ref{Global} is equivalent to the following condition
(refer to \cite{Kla86} for instance):
\begin{itemize}
\item[(a')] There exist some constants $A_{i, jk}^{\alpha\beta}$ and $B_{i, jk}^{ab, \alpha\beta}$ such that
\begin{align}
F_i^{\rm (W)} (w, \pa w, \pa_x\pa w)=& \sum_{\substack{1\le j, k\le N_2\\ 0\le |\alpha|, |\beta|\le 1}}
A_{i,jk}^{\alpha\beta} Q_0(\pa^\alpha w_j, \pa^\beta w_k)
\nonumber\\
& {}+\sum_{\substack{1\le j, k\le N_2\\ 0\le |\alpha|, |\beta|\le 1}}\sum_{0\le a<b\le3}
B_{i,jk}^{ab, \alpha\beta} Q_{ab}(\pa^\alpha w_j, \pa^\beta w_k)
\label{NullFormNullCond}
\end{align}
holds for any $i=N_1+1, \ldots, N$, and any $C^2$-function $w=(w_1,\ldots, w_{N_2})$,
where $\pa^\alpha=\pa_0^{\alpha_0}\pa_1^{\alpha_1}\pa_2^{\alpha_2}\pa_3^{\alpha_3}$
for a multi-index $\alpha=(\alpha_0, \alpha_1, \alpha_2, \alpha_3)$.
\end{itemize}

The condition (b) in Theorem \ref{Global} is assumed in order to compensate for bad behavior
of the solutions to the wave equations, as compared with their derivatives: 
Let $u=(v,w)$ be the solution to \eqref{OurSystem}.
The condition (b--i) says that if $k\in {\mathcal I}_1$, then all of $F_i^{\rm (q)}(u, \pa u, \pa_x \pa u)$ can
depend on $(\pa w_k, \pa_x\pa w_k)$, but {\it not} on $w_k$ itself
(remember that $\zeta_k(=\xi_{N_1+k})$ is the variable corresponding to $w_k(=u_{N_1+k})$),
while the divergence structure in the condition (b--ii)
assures that for each $k\in {\mathcal I}_2$, 
$w_k$ behaves better than we can expect in general 
(see Lemmas \ref{DivEne} and \ref{DivDecay} below; 
observe that the equation for $w_k(=u_{N_1+k})$ with $k\in {\mathcal I}_2$
in \eqref{OurSystem} is
$\dal w_k=F_{N_1+k}(u,\pa u, \pa_x\pa u)=\sum_{a=0}^3\pa_a\left\{{\mathcal G}_{k,a}(u, \pa u)\right\}$).
Here we remark that the condition (b--i) does not imply \eqref{b-ii-2} for $(l, k)\in {\mathcal I}_1\times {\mathcal I}_2$ in general, because we have
\begin{equation}
\label{TwoTen}
2Q_{ab}(\varphi, \psi)=\pa_a\left\{\varphi(\pa_b \psi)-(\pa_b\varphi)\psi\right\}
{}+\pa_b\left\{(\pa_a\varphi)\psi-\varphi(\pa_a \psi)\right\}
\end{equation}
for example (observe that $Q_{ab}(\varphi,\psi)$ depends only on $\pa \varphi$
and $\pa \psi$, while $\varphi(\pa_b\psi)-(\pa_b\varphi)\psi$ on the right-hand side
depends not only on $\pa\varphi$ and $\pa\psi$, but also on $\phi$ and $\psi$).

To help the understanding of our condition, we give a typical example here.
In what follows, for a finite family of functions 
$\{\phi_\lambda\}_{\lambda\in \Lambda}$ and a function $\psi$, we write 
$\psi=\sum_{\lambda\in \Lambda}' \phi_\lambda$
if there exists a family of
constants $\{c_\lambda\}_{\lambda\in \Lambda}$ such that $\psi=\sum_{\lambda\in \Lambda} c_\lambda \phi_\lambda$.
Let $u=(v,w)=(v, w_1, w_2)$ be an $\R^3$-valued function, 
and let $m$ be a positive constant.
Then the assumption in Theorem~\ref{Global} is fulfilled
with ${\mathcal I}_1=\{1\}$ and ${\mathcal I}_2=\{2\}$ 
for the following semilinear system: 
\begin{align}
(\dal+m^2) v = & \Sum_{|\alpha|, |\beta|\le 1}(\pa^\alpha v)(\pa^\beta v)
                      {}+\Sum_{\substack{|\alpha|\le 1 \\ 0\le b\le 3}} 
                      (\pa^\alpha v)(\pa_b w_1)+\Sum_{|\alpha|,|\beta|\le 1} (\pa^\alpha v)(\pa^\beta w_2)  \nonumber\\
                 & {}+\Sum_{0\le a,b\le 3} (\pa_a w_1)(\pa_b w_1)+\Sum_{\substack{0\le a\le 3\\|\beta|\le 1}}
                  (\pa_a w_1)(\pa^\beta w_2) \nonumber\\
                 &{}+\Sum_{|\alpha|,|\beta|\le 1} (\pa^\alpha w_2)(\pa^\beta w_2)+H_1(u, \pa u),
\\
\dal w_1 = & \Sum_{|\alpha|, |\beta| \le 1}(\pa^\alpha v)(\pa^\beta v)
                      {}+\Sum_{\substack{|\alpha|\le 1 \\ 0\le b\le 3}} 
                      (\pa^\alpha v)(\pa_b w_1)+\Sum_{|\alpha|,|\beta|\le 1} (\pa^\alpha v)(\pa^\beta w_2) \nonumber\\
& {}+\Sum_{j, k=1,2} Q_0(w_j, w_k){}+\Sum_{\substack{j,k=1,2 \\ 0\le a<b\le 3}} Q_{ab} (w_j, w_k) 
{}+H_2(u, \pa u), \\
\dal w_2 = & 
\sum_{a=0}^3 \pa_a \left(C_{1,a} v^2+ C_{2,a} v w_2 +H_{3,a}(u) \right)
{}+\Sum_{0\le a<b\le 3} Q_{ab}(v, w_2),
\label{Hinana2}
\end{align}
where $C_{j, a}$'s are real constants,  while
$H_1$, $H_2$, and $H_{3,a}$ are smooth functions in their arguments satisfying
$H_1(u, \pa u), H_2(u, \pa u)=O(|u|^3+|\pa u|^3)$ near $(u,\pa u)=0$, 
and $H_{3,a}(u)=O(|u|^3)$ near $u=0$. We use \eqref{TwoTen} to treat $Q_{ab}$ in \eqref{Hinana2}.

Now we would like to see the relation between our theorem and the previous results.
When $m_i>0$ for all $i=1, \ldots, N$, by regarding $v=u$, 
and by
neglecting the meaningless conditions (a) and (b),
Theorem \ref{Global} covers the previous results in \cite{Kla85:02} and \cite{Sha85} for nonlinear Klein-Gordon equations.
Similarly, when $m_i=0$ for all $i=1, \ldots, N$, 
by regarding $w=u$ and $N_1=0$ (thus $F_i^{\rm (W)}$ is regarded as $F_i^{\rm (q)}$ for $1\le i\le N$), it also covers the previous results in \cite{Chr86} and \cite{Kla86}  for nonlinear wave equations; note that the condition (b) for this case is automatically satisfied under the condition (a), 
because \eqref{NullFormNullCond} implies (b) with 
the choice of ${\mathcal I}_1=\{1,\ldots, N\}$
and ${\mathcal I}_2=\emptyset$.
It is easy to 
show that the strong null condition 
is satisfied if and only if
each $F_i^{\rm (q)}$ ($1\le i\le N$)
is a linear combination of $Q_{ab}(\pa^\alpha u_j, \pa^\beta u_k)$ with
$1\le j, k\le N$, $|\alpha|, |\beta| \le 1$ 
and $0\le a<b\le 3$.
Hence our conditions (a) and (b) are much weaker than the strong null condition
in \cite{Geo90}.
Note that some case of variable coefficients is also treated in \cite{Geo90},
but we can easily modify our conditions (a) and (b) to treat such case.

The main difficulty in the proof of Theorem \ref{Global} lies in the
fact that we can only use the vector fields which are 
compatible with both wave and Klein-Gordon equations.
To prove Theorem \ref{Global},
instead of the weighted $L^2$--$L^\infty$ estimate derived in \cite{Geo90},
we use weighted $L^\infty$--$L^\infty$ estimates
for wave equations (see Lemma~\ref{DecaySolWave} below), which require
a smaller set of vector fields than the admissible set of vector fields for
the Klein-Gordon equations. 
We also need some estimates for null forms without using the scaling operator $S$;
they will be given in Lemma~\ref{SKnull} below.
To treat $F_i^{\rm (W)}$ with $1\le i\le N_1$, for which the null condition is {\it not} assumed,
we adopt a technique used in Y.~Tsutsumi~\cite{YTsu03}, where
the Dirac-Proca equations are considered (see \eqref{DP-D}--\eqref{DP-P} below).
This technique is motivated by Bachelot~\cite{Bac88} and Kosecki~\cite{Kos92}, and it is
closely related to the normal form technique in Shatah~\cite{Sha85}. 
We will prove Theorem \ref{Global} in Section~5.
 
We conclude this section with some examples from physics which can be treated by
Theorem \ref{Global}. Note that all the following examples are semilinear
(or can be regarded as semilinear), and the conditions \eqref{quasi-linear} and
\eqref{symmetricity} in Theorem \ref{Global} are trivially satisfied.
Thus we only have to check the conditions (a) and (b). 
%The initial data are always supposed to be sufficiently small.
%%%%%%%%%%%%
\begin{example}[The Dirac-Klein-Gordon equations]
\normalfont
Let us consider the Dirac equation coupled with the Klein-Gordon
or wave equation:
\begin{align}
\label{DKG-D}
-\sqrt{-1}\sum_{a=0}^3 \gamma_a\pa_a\psi+M\psi=&\sqrt{-1}c\varphi \gamma_5 \psi,\\
\left(\dal+m^2\right)\varphi=& \psi^* H \psi
\label{DKG-KG}
\end{align}
in $(0,\infty)\times \R^3$, where $\sqrt{-1}$ denotes the imaginary unit, $M, m\ge 0$, $c$ is a real constant, $H$ is a $4\times 4$ Hermitian matrix, 
$\psi$ is a $\C^4$-valued function, $\varphi$ is a real valued function,
and $\psi^*$ denotes the complex conjugate transpose of $\psi$.
$\gamma_a$ $(0\le a\le 3$) in the above are $4\times 4$ matrices satisfying 
$\gamma_a\gamma_b+\gamma_b\gamma_a=2g_{ab}I$ for $0\le a,b\le 3$,
where $I$ is the $4\times 4$ identity matrix, and $(g_{ab})_{0\le a, b\le 3}=\diag (1,-1,-1,-1)$;
$\gamma_5$ is defined by $\gamma_5=-\sqrt{-1}\gamma_0\gamma_1\gamma_2\gamma_3$.
The initial data are supposed to be sufficiently small.

We set $D_M^{\pm}=\pm \sqrt{-1}\sum_{a=0}^3\gamma_a\pa_a+MI$. Since we have $\gamma_a\gamma_5=-\gamma_5\gamma_a$ for $0\le a\le 3$, we get
$D_M^{+}\gamma_5=\gamma_5D_M^{-}$. Therefore, operating $D_M^{+}$ to \eqref{DKG-D}, we get
\begin{align}
\left(\dal+M^2\right) \psi=& \sqrt{-1}cD_M^{+}(\varphi \gamma_5\psi)=-c\sum_{a=0}^3
(\pa_a\varphi)\gamma_a\gamma_5\psi+\sqrt{-1}c\varphi\gamma_5 
(D_M^{-}\psi) \nonumber\\
=& -c\sum_{a=0}^3(\pa_a\varphi)\gamma_a\gamma_5\psi-c^2\varphi^2\psi,
\label{DKG-D-Red}
\end{align}
where we have used \eqref{DKG-D} and $\gamma_5\gamma_5=I$ to obtain the last identity.
If $M>0$ and $m>0$, the system \eqref{DKG-KG}--\eqref{DKG-D-Red} is a system
of the nonlinear Klein-Gordon equations, and we have the global solution.
When $M>0$ and $m=0$, putting $u=(v,w)$ with $v=(\RealPart \psi, \ImPart \psi)$ and $w(=w_1)=\varphi$,
we see that 
the conditions (a) and (b) 
are satisfied for the system \eqref{DKG-KG}--\eqref{DKG-D-Red},
and thus Theorem \ref{Global} is applicable; more precisely the condition (b) is satisfied with ${\mathcal I}_1=\{1\}$ and ${\mathcal I}_2=\emptyset$.
The global existence result for this case where $M>0$ and $m=0$, with compactly supported
initial data, has been already obtained
by Bachelot \cite{Bac88}.
Differently from \cite{Bac88}, we can also treat the case where $M=0$ and $m>0$. Indeed, 
the first identity in \eqref{DKG-D-Red} with $M=0$ can be read as
$\dal \psi=-c\sum_{a=0}^3 \gamma_a \pa_a (\varphi \gamma_5 \psi)$, and putting $u=(u_i)_{1\le i\le 9}=(v, w)$
with $v(=v_1)=\varphi$ and $w=(w_k)_{1\le k\le 8}=(\RealPart \psi, \ImPart \psi)$, we can verify the conditions (a)
and (b) with ${\mathcal I}_1=\emptyset$ and ${\mathcal I}_2=\{1, \ldots, 8\}$.
This last case is closely connected to the next example, the Dirac--Proca equations.
\end{example}
%%%%%%%%%%%%%%%%%%%%%%%%%%
\begin{example}[The Dirac-Proca equations]
\normalfont
Y.~Tsutsumi~\cite{YTsu03} proved the global existence of small solutions to
the Dirac-Proca equations, which can be reduced to the following
coupled system of the massless Dirac and the Klein-Gordon equations:
\begin{align}
\label{DP-D}
-\sqrt{-1} \sum_{a=0}^3 \gamma_a\pa_a \psi=& -\frac{1}{2}
\sum_{a=0}^3 g_{aa} A_a %\varphi_a
                        \gamma_a(I+\gamma_5)\psi,\\
\label{DP-P}
(\dal+m^2) A_a 
 =&\frac{1}{2} \psi^* \gamma_0 \gamma_a (I+\gamma_5)\psi,
\quad a=0,1,2,3
\end{align}
with the constraint $\sum_{a=0}^3 \pa_a A_a %\varphi_a
                                                        =0$ at $t=0$,
where $m>0$, $\psi$ is a $\C^4$-valued function, 
and $A_a %\varphi_a
$ for $0\le a\le 3$
are real-valued functions.
%%%%%%%%%%%%%%%%%%%%%%%% 
$I$, $\gamma_a$ ($a=0,1,2,3,5)$ and $(g_{ab})$ are as in
the Dirac-Klein-Gordon equations.
In a similar manner to \eqref{DKG-D-Red} with $M=0$, from \eqref{DP-D} we obtain
\begin{equation}
\label{DP-D-Red}
\dal \psi=-\frac{1}{2}\sqrt{-1} \sum_{b=0}^3 \pa_b \left(
\sum_{a=0}^3 g_{aa} A_a %\varphi_a
     \gamma_b\gamma_a(I+\gamma_5)\psi
\right).
\end{equation}
Putting $u=(v,w)\in \R^4\times \R^8$ with $v=(A_a %\varphi_a
)_{0\le a\le 3}$, and $w=(\RealPart \psi, \ImPart \psi)$,
we find that the conditions (a) and (b) 
hold for the system \eqref{DP-P}--\eqref{DP-D-Red},
and thus Theorem \ref{Global} is applicable;
more precisely \eqref{DP-D-Red} implies that (b) is satisfied 
with ${\mathcal I}_1=\emptyset$ and ${\mathcal I}_2=\{1, \ldots, 8\}$.
\end{example}
%%%%%%%%%%%%%%%%%%%%%%%%%%%
\begin{example}[The Klein-Gordon-Zakharov equations]
\normalfont
Ozawa-Tsutaya\break
-Tsutsumi~\cite{OzaTsuTsu95} and Tsutaya~\cite{KTsu96}
proved the global existence of small solutions to
the Klein-Gordon-Zakharov equations:
\begin{equation}
\label{KGZ}
\begin{cases}
(\dal+1) \widetilde{u}= -\widetilde{n}\widetilde{u},\\
\dal \widetilde{n}=  \Delta_x |\widetilde{u}|^2
\end{cases}
\end{equation}
in $(0,\infty)\times \R^3$,
where $\widetilde{u}=(\widetilde{u}_1,\widetilde{u}_2, \widetilde{u}_3)$ is a $\C^3$-valued function, and $\widetilde{n}$ is a real valued function
(see also Ozawa-Tsutaya-Tsutsumi~\cite{OzaTsuTsu99} for the multiple speed case).
 
By setting
$v_i=\RealPart \widetilde{u}_i$, $v_{i+3}=\ImPart \widetilde{u}_i$ ($1\le i\le 3$),
$v_{3i+3+k}=\pa_k v_i $ ($1\le i\le 6,\ 1\le k\le 3$),
and $w(=w_1)=\widetilde{n}$,
we see that solving \eqref{KGZ} is equivalent 
to solving
\begin{equation}
\label{ReducedKGZ}
\begin{cases}
(\dal+1)v_i=-wv_i, \qquad\qquad\qquad\qquad\quad\  \,
1\le i\le 6,\\
(\dal+1) v_{3i+3+k}
=-w(\pa_k v_i)-(\pa_k w)v_i,
\quad 1\le i\le 6, 1\le k\le 3,\\
\dal w
=\sum_{j=1}^3 \pa_j \sum_{i=1}^3 2(v_i v_{3i+3+j}+v_{i+3}v_{3i+12+j}).
\end{cases}
\end{equation}
Note that the system \eqref{ReducedKGZ} is a semilinear system
of $$
u=(u_1, \ldots, u_{25})=(v_1,\ldots, v_{24}, w_1)=(v,w).
$$
The conditions (a) and (b) (with ${\mathcal I}_1=\emptyset$ and ${\mathcal I}_2=\{1\}$)
are satisfied for \eqref{ReducedKGZ}.
% follows from $F^W\equiv 0$,
%  and we can verify the condition (b)
% by putting 
% ${\mathcal I}_1=\emptyset$, ${\mathcal I}_2=\{1\}$,
% $$
% {\mathcal G}_{1, j}\bigl((\eta, \zeta), (\eta', \zeta')\bigr)=2\sum_{i=1}^3 (\eta_i \eta_{3i+3+j}+\eta_{i+3}\eta_{3i+12+j})
% $$
% for $j=1,2,3$, and ${\mathcal G}_{1,0}=0$.
Hence
we can apply Theorem \ref{Global} to
show the global existence of small amplitude
solutions to \eqref{KGZ}.
%%%%%%
\end{example}
%%%%%%%%%%%%%%%%%%%%%%%%%%%%%%%%%%%%%%%%%%%%%%%%%
\begin{example}
\normalfont
The last example is not from physics,
as far as the author knows. 
This example shows that some change of unknowns may help us to apply our theorem.
Consider
\begin{equation}
\label{KataEx}
\begin{cases}
(\dal+1)v= w^2,\\
\dal w= v^2
\end{cases}
\end{equation}
in $(0,\infty)\times \R^3$.
We can treat this example in the following way,
though it does not explicitly satisfy the assumption
of Theorem \ref{Global}:
Set $\widetilde{v}=v-w^2$ (cf.~\eqref{YTran} below). Then 
we get
\begin{equation}
\label{KataExR}
\begin{cases}
(\dal+1) \widetilde{v}= -2Q_0(w, w)
{}-2w\left(\widetilde{v}+w^2\right)^2,\\
\dal w=\left(\widetilde{v}+w^2\right)^2
\end{cases}
\end{equation}
(cf.~\eqref{YTran03} below).
This system \eqref{KataExR}
satisfies the assumption in Theorem~\ref{Global} with $u=(\widetilde{v}, w)=(\widetilde{v}_1, w_1)$, $I_1=\{1\}$, and $I_2=\emptyset$.
Thus we get a global solution $(\widetilde{v}, w)$ to \eqref{KataExR} for small data, and accordingly we obtain a
global solution $(v, w)$ to the original system \eqref{KataEx}. 
\end{example}
%%%%%%%%%%%%%%%%%%%%%%%%%%%%%%%%%%%%%%%%%%%
\section{Preliminary Results}
%%%%%%%%%%%%%%%%%%%%%%%%%%%%%%%%%%%%%%%%%%%
In this section, we state the known estimates
for the wave and Klein-Gordon equations.
Throughout this paper,
we write $\jb{z}=\sqrt{1+|z|^2}$ for $z\in \R^d$, where $d$ is a positive integer.

We start this section with the energy inequality for hyperbolic systems,
which can be shown easily by the standard method.
%%%%%%%%%%%%%%%%%%%%%%%%%%%%%
\begin{lemma}\label{EnergyHyp}
 Let $m_i\ge 0$ for $1\le i\le N$, and $T>0$.
Suppose that 
 $\widetilde{\gamma}=\left(\widetilde{\gamma}^{\, ij}_{ka}\right)$
be a smooth function satisfying
$$
\widetilde{\gamma}^{\, ij}_{ka}(t,x)=\widetilde{\gamma}_{ka}^{\, ji}(t,x),
\ \widetilde{\gamma}^{\, ij}_{kl}(t,x)=\widetilde{\gamma}^{\, ij}_{l k}(t,x),\ (t,x)\in (0, T)\times \R^3
$$
for $1\le i, j\le N$, $1\le k, l \le 3$,  and $0\le a\le 3$. 
We also assume that
$$
\left|\sum_{1\le i, j\le N} \sum_{1\le k,l\le 3} \widetilde{\gamma}_{kl}^{\, ij}(t,x) \xi'_{i,k} \xi'_{j,l}\right|
\le \frac{1}{2} \sum_{1\le i \le N}\sum_{1\le k\le 3} |\xi'_{i,k}|^2,\ (t,x)\in (0, T)\times \R^3
$$
for any $%\widetilde{\xi}'=
(\xi'_{i,k})_{1\le i\le N, 1\le k\le 3}\in \R^{3N}$.

Let $\varphi=(\varphi_1, \ldots, \varphi_N)$ be the solution to
\begin{align*}
& \left(
 \dal+m_i^2\right) \varphi_i-\sum_{1\le j\le N} \sum_{\substack{1\le k\le 3\\ 0\le a\le 3}} \widetilde{\gamma}_{ka}^{\, ij}(\pa_k \pa_a \varphi_j)=\Phi_i \quad  \text{ in } (0, T)\times \R^3,\\
& \varphi_i(0,x)=\varphi_{i}^{(0)}(x), \ (\pa_t \varphi_i)(0,x)=\varphi_{i}^{(1)}(x), \quad x\in \R^3
\end{align*}
for $1\le i\le N$,
where 
$\varphi^{(0)}=(\varphi^{(0)}_1, \ldots, \varphi^{(0)}_N) %_{1\le i\le N}
\in H^1(\R^3; \R^N)$, 
$\varphi^{(1)}=(\varphi_{1}^{(1)},\ldots, \varphi^{(1)}_N) %_{1\le i\le N} 
  \in L^2(\R^3; \R^N)$,
and
$\Phi=(\Phi_1, \ldots, \Phi_N) \in L^1\bigl((0,T); L^2(\R^3; \R^N)\bigr)$.
Then, there exists a positive constant $C$, which is independent of $T$, such that
\begin{align*}
& \sum_{i=1}^N \left(\|\pa \varphi_i(t)\|_{L^2}+m_i\|\varphi_i(t)\|_{L^2}\right) \\
&  \quad \le C\left(\bigl\|\varphi^{(0)}\bigr\|_{H^1}+\bigl\|\varphi^{(1)}\bigr\|_{L^2}
{}+\int_0^t \left\|\pa\widetilde{\gamma}(\tau)\right\|_{L^\infty}\|\pa \varphi(\tau)\|_{L^2}d\tau
{}+\int_0^t \|\Phi(\tau)\|_{L^2} d\tau\right)
\nonumber
\end{align*}
for $0\le t<T$.
\end{lemma}
%%%%%%%%%%%%%%%%%%%%%%%%%%%%%%%%%%%%%%

Before we proceed to the decay estimates of the solutions
to the Klein-Gordon and wave equations,  we introduce the vector fields $\Omega_j$ and $L_j$
for $1\le j\le 3$ by
\begin{align}
\label{frame01}
\Omega=& (\Omega_1, \Omega_2, \Omega_3)=x\times \nabla_x=(x_2\pa_3-x_3\pa_2, x_3\pa_1-x_1\pa_3,
x_1\pa_2-x_2\pa_1),\\ 
\label{frame01'}
L=& (L_1, L_2, L_3)=x\pa_t+t\nabla_x=(x_1\pa_t+t\pa_1, x_2\pa_t+t \pa_2, x_3\pa_t+t\pa_3),
\end{align}
where $\nabla_x=(\pa_1, \pa_2, \pa_3)$, and $\times$ is the external product in $\R^3$.
Writing $\pa=(\pa_a)_{0\le a\le 3}$, we set
\begin{equation}
\label{KGVec}
Z=(Z_1, \ldots, Z_{10})=\bigl((\Omega_j)_{1\le j\le 3}, (L_j)_{1\le j\le 3}, (\pa_a)_{0\le a\le 3}\bigr)
=(\Omega, L, \pa).
\end{equation}
Note that we have
\begin{equation}
\left[L_j, \dal+m^2\right]
=\left[\Omega_j, \dal+m^2\right]
=\left[\pa_a, \dal+m^2\right]=0,
\ 1\le j\le 3,\ 0\le a\le 3
\end{equation}
for $m\ge 0$, where $[A,B]=AB-BA$ for operators $A$ and $B$. 
Hence the vector fields in $Z=(\Omega, L, \pa)$ are compatible with the Klein-Gordon equations,
as well as the wave equations.
Here we note that we have
\begin{equation}
\label{CommZD}
[Z_j, \pa_a]=\sum_{b=0}^3 C^{ja}_b \pa_b,\ [Z_j, Z_k]=\sum_{l=1}^{10} D_l^{jk} Z_l,
\quad 1\le j, k\le 10,\ 0\le a\le 3
\end{equation}
with appropriate constants $C^{ja}_b$ and $D_l^{jk}$.
For a multi-index $\alpha=(\alpha_1, \ldots, \alpha_{10})$,
we write $Z^\alpha=Z_1^{\alpha_1}\cdots Z_{10}^{\alpha_{10}}$.
For a function $\varphi=\varphi(t,x)$ and a nonnegative integer $s$,
we define
\begin{equation}
\label{InvariantNorm}
|\varphi(t,x)|_s=\sum_{|\alpha|\le s} |Z^\alpha \varphi(t,x)|,
\quad \|\varphi(t)\|_s=\sum_{|\alpha|\le s} \left\|Z^\alpha \varphi(t,\cdot)\right\|_{L^2(\R^3)}.
\end{equation}

Using the vector fields in $Z=(\Omega, L, \pa)$, Klainerman~\cite{Kla85:02} obtained the decay estimate
for the solutions to the Klein-Gordon equations. This estimate has been modified and generalized
by many authors (for instance, see Bachelot~\cite{Bac88}, H\"ormander~\cite{Hoe97}, Sideris~\cite{Sid89}, and Georgiev~\cite{Geo92}).
Here we state the estimate obtained in \cite{Geo92}:
Let $\chi_j$ ($j\ge 0$) be nonnegative $C^\infty_0(\R)$-functions 
satisfying
\begin{align}
& \sum_{j=0}^\infty \chi_j(\tau)=1\text{ for $\tau\ge 0$,}\\
& \text{$\supp \chi_j=[2^{j-1}, 2^{j+1}]$ for $j\ge 1$, and $\supp \chi_0\cap [0,\infty)=[0,2]$.}
\end{align}
%%%%%%%%%%%%%%%%%%%%%%%%%%%%%%%%%%%%%%%
\begin{lemma}\label{DecayKleinGordon}
Let $m>0$, and $v$ be a smooth solution to
$$
\left(\dal+m^2\right)v(t,x)=\Phi(t,x),\ (t,x)\in (0,\infty)\times \R^3.
$$
Then there exists a positive constant $C=C(m)$ such that we have
\begin{align}
\jb{t+|x|}^{3/2}|v(t,x)|
\le & C \sum_{j=0}^\infty \sum_{|\alpha|\le 4} \sup_{\tau\in [0,t]} \chi_j(\tau) 
\left\|\jb{\tau+|\cdot|} Z^\alpha \Phi(\tau,\cdot)\right\|_{L^2(\R^3)}
\nonumber\\
&  
{}+C \sum_{j=0}^\infty \sum_{|\alpha|\le 5}
\left\|\jb{\,\cdot\,}^{3/2}\chi_j(|\cdot|)Z^\alpha v(0,\cdot)\right\|_{L^2(\R^3)}
\label{GeoEs}
\end{align}
for $(t,x)\in (0, \infty)\times \R^3$,
provided that the right-hand side of \eqref{GeoEs} is finite. 
\end{lemma}
%%%%%%%%%%%%%%%%%%%%%%%%%%%%%%%%%%%%%%%%%%%%
For the proof, see Georgiev~\cite[Theorem 1]{Geo92}.
%%%%%%%%%%%%%%%%%%%%%%%%%%%%%%%%%%%%%%%%%%%%
\smallskip

Now we turn our attention to the wave equations.
In \cite{Kla86}, a weighted $L^1$--$L^\infty$ estimate 
for the wave equation
is derived (see also H\"ormander \cite{Hoe88}), where
the scaling operator $S=t\pa_t+x\cdot \nabla_x$ as well as $Z=(\Omega, L, \pa)$
is used.
Since we have $[S, \dal]=-2\dal$, the scaling operator 
$S$ is applicable to the wave equations,
but it is incompatible with the Klein-Gordon equations.
Therefore Georgiev (\cite{Geo90}) developed a weighted $L^2$--$L^\infty$ estimate
involving only $Z$.
%On the other hand, t
There is also a large literature on the study of systems
of nonlinear wave equations with multiple speeds of the form
$$
 \dal_{c_i} u_i=F_i(u, \pa u, \pa_x \pa u),\quad 1\le i\le N,
$$
where $c_i>0$ ($1\le i\le N$) and $\dal_c=\pa_t^2-c^2\Delta_x$
(see, for example, \cite{Kat-Kub08:01} and the references cited therein).
In the study of this kind of system, the vector field method
using only $(S, \Omega, \pa)$ has been developed, 
because $L=(L_j)_{1\le j\le 3}$ is incompatible with such system 
(observe that $[L_j, \dal_c]=2(c^2-1)\pa_t\pa_j$ has no good property when $c\ne 1$).
Especially, in Yokoyama~\cite{Yok00} and Kubota-Yokoyama~\cite{Kub-Yok01} (see also the author~\cite{Kat04:02}),
weighted $L^\infty$--$L^\infty$ estimates requiring only $(\Omega, \pa)$
are adopted to prove some global existence results under the null condition
(the origin of these estimates can be found 
in John~\cite{Joh79} and Kovalyov~\cite{Kov87}; see also Kovalyov-Tsutaya~\cite{Kov-Tsu}).
We will employ these $L^\infty$--$L^\infty$ estimates in the proof of Theorem~\ref{Global}
because they require only $(\Omega, \pa)$ and are easily applicable 
to the coupled system of the wave and Klein-Gordon equations.
Here we note that $S$ is still used in the arguments in \cite{Yok00}, \cite{Kub-Yok01} and \cite{Kat04:02}
to treat the null forms (see \eqref{YokNull} below).

To state the weighted $L^\infty$--$L^\infty$ estimates, 
we define
\begin{align}
\label{KY-weight}
{\mathcal W}_\rho(t, r):=&
 \begin{cases}
  \jb{t+r}^{\rho} & \text{if $\rho<0$},\\
  \left\{\log \left(2+\jb{t+r}\jb{t-r}^{-1}\right)
  \right\}^{-1} & \text{if $\rho=0$},\\
  \jb{t-r}^{\rho} & \text{if $\rho>0$}.                  
 \end{cases}
\end{align}
We also introduce
\begin{equation}
\label{WeakWei}
W_-(t,r):=\min\left\{\jb{r}, \jb{t-r} \right\}.
\end{equation}

For the homogeneous wave equations we have the following estimate
which was essentially proved in Asakura~\cite[Proposition 1.1]{Asa86} (see also \cite[Lemma 3.1]{Kat-Yok06} for the expression below):
\begin{lemma}\label{Asakura}
Let $w$ be a smooth solution to
$$
\dal w(t,x)=0,\ (t,x)\in (0,\infty)\times \R^3
$$
with initial data $w=w^{(0)}$, $\pa_t w=w^{(1)}$ at $t=0$.

Let $\kappa>0$.
Then, there exists a positive constant $C=C(\kappa)$ such that
\begin{align}
& 
\jb{t+|x|} 
{\mathcal W}_{\kappa-1}(t,|x|) |w(t,x)|
\nonumber\\
& \qquad \le
C \sup_{|y-x|\le t} \jb{y}^\kappa \left(\jb{y} \sum_{|\alpha|\le 1} \bigl|\bigl(\pa_x^\alpha w^{(0)}\bigr)(y)\bigr|
{}+|y| \bigl|w^{(1)}(y)\bigr|\right)
\label{Asakura01}
\end{align}
for $(t,x)\in (0,\infty)\times \R^3$.
Here $\pa_x=(\pa_1, \pa_2, \pa_3)$, and we have used the standard notation of
multi-indices.
\end{lemma}

The following weighted $L^\infty$--$L^\infty$ estimates 
are the special cases of the estimates obtained in Kubota-Yokoyama~\cite[Lemma 3.2]{Kub-Yok01}
(see also Katayama-Kubo~\cite[Lemma 3.4]{Kat-Kub08:03} for the expression below):
%%%%%%%%%%%%%%%%%%%%%%%%%%%%%%%%%%%%%%%%%%%%%
\begin{lemma}\label{DecaySolWave}
Let $w$ be a smooth solution to
$$
\dal w(t,x)=\Psi(t,x),\ (t,x)\in (0,\infty)\times \R^3
$$
with initial data $w=\pa_t w=0$ at $t=0$.

Suppose that $\rho\ge 0$, $\kappa\ge 1$, and $\mu>0$.
Then there exists a positive constant $C=C(\rho, \kappa, \mu)$ such that
\begin{align}
& \jb{t+|x|}^{1-\rho}{\mathcal W}_{\kappa-1}(t,|x|)
|w(t,x)| \nonumber\\
& \quad \le  
C\sup_{\tau\in [0,t]}\sup_{|y-x|\le t-\tau} |y|\jb{\tau+|y|}^{\kappa-\rho+\mu}
W_-(\tau, |y|)^{1-\mu}|\Psi(\tau,y)|,          
\label{KubYok}\\
& \jb{t+|x|}^{-\rho} \jb{x} \jb{t-|x|}^{\kappa}
|\pa w(t,x)| \nonumber\\
& \quad \le  
C\sup_{\tau\in [0,t]}\sup_{|y-x|\le t-\tau} |y|\jb{\tau+|y|}^{\kappa-\rho+\mu} W_-(\tau, |y|)^{1-\mu}
\sum_{|\alpha|+|\beta|\le1}|\pa^\alpha\Omega^\beta \Psi(\tau,y)|          
\label{KubYokD}
\end{align}
for $(t,x)\in (0,\infty)\times \R^3$.
Here $\pa=(\pa_0, \pa_1, \pa_2, \pa_3)$, and $\Omega$ is given by \eqref{frame01}.
\end{lemma}
%%%%%%%%%%%%%%%%%%%%
The following Sobolev type inequality will be used to combine decay estimates with the energy estimates
(see Klainerman~\cite{Kla87} for the proof):
%%%%%%%%%%%%%%%%%%%%%%%%%%%%%%%%%%%%%%%%%%%%
\begin{lemma}\label{KlaSobolev}
For a smooth function $\varphi$ on $\R^3$, we have
\begin{equation}
\label{KlaSov}
\sup_{x\in \R^3} \jb{x} |\varphi(x)|\le C \sum_{|\alpha|+|\beta|\le 2} 
\left\|\pa_x^\alpha \Omega^\beta \varphi\right\|_{L^2(\R^3)},
\end{equation}
provided that the right-hand side of \eqref{KlaSov} is finite.
Here $C$ is a universal positive constant.
\end{lemma}
%%%%%%%%%%%%%%%%%%%%%%%%%%%%%%%%%%%%%%%%%%%

We conclude this section with some observation on the wave equations of the following type:
\begin{equation}
\label{DivWav}
\begin{cases}
\dal \psi(t,x)=\sum_{a=0}^3 \pa_a \Psi_a(t,x), & (t,x)\in (0,\infty)\times \R^3,\\
\psi(0,x)=\psi^{(0)}(x),\ (\pa_t \psi)(0,x)=\psi^{(1)}(x), & x\in \R^3.
\end{cases}
\end{equation}
For $0\le a\le 3$, let $\psi_a=\psi_a(t,x)$ be the solution to
$\dal \psi_a=\Psi_a$ with initial data $\psi_a=\pa_t \psi_a=0$ at $t=0$, 
and let $\psi_{\rm f}(t,x)$ be the solution to $\dal \psi_{\rm f}=0$
with initial data $\psi_{\rm f}=\psi^{(0)}$ 
and $\left(\pa_t \psi_{\rm f}\right)=\psi^{(1)}-\Psi_0(0,\cdot)$ at $t=0$.
It is easy to verify that the solution $\psi$ to \eqref{DivWav} can be written as
$\psi=\sum_{a=0}^3 \pa_a \psi_a+\psi_{\rm f}$.
Therefore, we can essentially regard $\psi$ as
derivatives of solutions to some wave equations, and $\psi$ enjoys better estimates than we can expect
in general.
%%%%%%%%%%%%%%%%%%%%%%%%%%%%
\begin{lemma}\label{DivEne}
Let $\psi$ be the solution to \eqref{DivWav}. Then we have
\begin{align}
 \|\psi(t,\cdot)\|_{L^2(\R^3)}
 \le &C\Bigl(\bigl\|\psi^{(0)}\bigr\|_{L^2(\R^3)}+\bigl\|\psi^{(1)}\bigr\|_{L^{6/5}(\R^3)}
{}+\left\|\Psi_0(0,\cdot)\right\|_{L^{6/5}(\R^3)}\Bigr)
\nonumber\\
     &{}+C\sum_{a=0}^3\int_0^t \left\|\Psi_a(\tau,\cdot)\right\|_{L^2(\R^3)} d\tau,
\label{DivEneS}
\end{align}
provided that the right-hand side of \eqref{DivEneS} is finite.
\end{lemma}
%%%%%%%%%%%%%%%%%%%%%%%%%%%%
{\it Proof.}
We have $\|\psi\|_{L^2}\le \sum_{a=0}^3\|\pa \psi_a\|_{L^2}+\|\psi_{\rm f}\|_{L^2}$.
Hence \eqref{DivEneS} follows from the energy inequality (cf.~Lemma \ref{EnergyHyp})
for $\psi_a$ ($0\le a\le 3$), and
the $L^2$-estimate for $\psi_{\rm f}$ (see Strauss \cite{Str68} for example).
%\qed
\begin{lemma}\label{DivDecay}
Let $\psi$ be the smooth solution to \eqref{DivWav}. 
Suppose that $\rho\ge 0$, $\kappa\ge 1$ and $\mu>0$. Then we have
\begin{align}
& \jb{t+|x|}^{-\rho}\jb{x}\jb{t-|x|}^{\kappa}|\psi(t,x)| \nonumber\\
& \quad
\le C \sup_{\tau\in [0,t]}\sup_{|y-x|\le t-\tau}
|y|\jb{\tau+|y|}^{\kappa-\rho+\mu}W_-(\tau,|y|)^{1-\mu}
\sum_{\substack{|\alpha|+|\beta|\le 1\\ 0\le a\le 3}}
|\pa^\alpha \Omega^\beta \Psi_a(\tau,y)|\nonumber\\
& \ \qquad {}+C\sup_{|y-x|\le t} \jb{y}^{\kappa+1-\rho}
\left(
 \jb{y}\sum_{|\alpha|\le 1} \bigl|\bigl(\pa_x^\alpha \psi^{(0)}\bigr)(y)\bigr|+|y|\,\bigl|\psi^{(1)}(y)\bigr| 
\right).
\label{KubYokO}
\end{align}
\end{lemma}
{\it Proof.}
The estimates for $\pa_a \psi_a$ 
follow from \eqref{KubYokD};
using Lemma~\ref{Asakura} to estimate 
$\psi_{\rm f}$, and noting that
we have 
$$
\jb{y}^{\kappa+1-\rho}|y|\,|\Psi_0(0,y)|=\left\{|y|\jb{\tau+|y|}^{\kappa-\rho+\mu}
 W_-(\tau,|y|)^{1-\mu} |\Psi_0(\tau,y)|\right\}\bigr|_{\tau=0},
$$ 
we obtain the desired result immediately
(more precisely, \eqref{KubYokO} for zero initial data is 
directly proved in \cite{Kub-Yok01}, and
\eqref{KubYokD} is obtained as its corollary in fact).
%\qed 
%\end{proof}

%%%%%%%%%%%%%%%%%%%%%%%%%%%%%%%%%%%%%%%%%%%
%%%%%% Null Forms
%%%%%%%%%%%%%%%%%%%%%%%%%%%%%%%%%%%%%%%%%%%
\section{Estimates for the Null Forms}\label{NullFormEstimate}
In this section, we will derive some estimates for the null forms.
We set
$r=|x|$,  $\omega=(\omega_1,\omega_2,\omega_3)$ with $\omega_j=x_j/r$ for $1\le j\le 3$,
and $\pa_r=\omega\cdot \nabla_x=\sum_{j=1}^3 \omega_j \pa_j$.
Then we have
\begin{equation}
\nabla_x= \omega \pa_r-r^{-1} 
(\omega\times \Omega),
\label{frame02}
\end{equation}
where 
$\Omega$ is defined by \eqref{frame01}. Since \eqref{frame01} and \eqref{frame01'} yield
$tr^{-1} \Omega=\omega \times L$,
the expression \eqref{frame02} implies
\begin{equation}
\label{frame04}
(t+r)(\nabla_x-\omega\pa_r)=-\omega\times(\Omega+\omega\times L).
\end{equation}
From \eqref{frame04}, we obtain
\begin{align}
&\left|Q_0(\varphi, \psi)-Q_0^{\rm rad}(\varphi, \psi)\right|
{}+\sum_{k=1}^3\left|Q_{0k}(\varphi, \psi)-\omega_k Q_{0r}^{\rm rad}(\varphi, \psi)\right|
\nonumber\\
&\qquad {}+\sum_{1\le j<k\le 3} \left|Q_{jk}(\varphi, \psi)\right|
\le C\jb{t+r}^{-1}\left(|Z \varphi|\,|\pa \psi|+|\pa \varphi|\, |Z\psi|\right) 
\label{NullRep01}
\end{align}
at $(t,x)\in [0,\infty)\times \R^3$,
where 
$$
Q_0^{\rm rad}(\varphi, \psi):=(\pa_t \varphi)(\pa_t \psi)- (\pa_r \varphi)(\pa_r \psi), 
\ Q_{0r}^{\rm rad}(\varphi, \psi):=(\pa_t \varphi)(\pa_r \psi)- (\pa_r \varphi)(\pa_t \psi),
$$
$Z$ is given by \eqref{KGVec}, and $C$ is a positive constant.
Putting $L_r:=\omega\cdot L=r\pa_t+t\pa_r$,
we get
\begin{align}
\label{NullRep02}
(t+r)Q_{0r}^{\rm rad}(\varphi, \psi)=
(\pa_t \varphi-\pa_r\varphi)(L_r\psi)-(L_r\varphi)(\pa_t\psi-\pa_r\psi).
\end{align}
From \eqref{NullRep01} and \eqref{NullRep02}, we obtain
\begin{equation}
\label{NullDecay01}
|Q_{ab}(\varphi, \psi)|\le
C\jb{t+r}^{-1}\left(|Z\varphi|\, |\pa \psi|+
|\pa\varphi|\, |Z\psi|\right)
\end{equation}
for $0\le a<b\le 3$
at $(t,x)\in [0,\infty)\times \R^3$,
where $C$ is some positive constant.
Thus we only need the vector fields in $Z=(\Omega, L, \pa)$ to obtain the extra decay
factor $\jb{t+r}^{-1}$ for terms satisfying the strong null
condition.

To treat the null form $Q_0$, we introduce $\pa_{\pm}:=\pa_t\pm \pa_r$.
Then we get
\begin{equation}
\label{NullRep03}
Q_0^{\rm rad}(\varphi, \psi)=\frac{1}{2}\left((\pa_+\varphi)(\pa_-\psi)+(\pa_- \varphi)(\pa_+\psi)\right).
\end{equation}
As we will see below, estimates of
$\pa_+\varphi$ and $\pa_+\psi$ are important in deriving enhanced decay for $Q_0$. 
Note that we also have
\begin{equation}
\label{NullRep04}
Q_{0r}^{\rm rad}(\varphi, \psi)=(\pa_+\varphi)(\pa_r\psi)-(\pa_r \varphi)(\pa_+\psi).
\end{equation}

Rewriting $\pa_+$ as 
\begin{equation}
\pa_+=(t+r)^{-1}(S+L_r)
\end{equation}
with $S=t\pa_t+x\cdot\nabla_x=t\pa_t+r\pa_r$,
from \eqref{NullRep01} and \eqref{NullRep03} we obtain
\begin{equation}
\label{NullDecay02}
|Q_0(\varphi, \psi)|\le C\jb{t+r}^{-1}
\left(|\Gamma \varphi|\,|\pa \psi|+|\pa \varphi|\,|\Gamma \psi|
\right),
\end{equation}
where $\Gamma:=(S,Z)=(S,L,\Omega, \pa)$.
The estimate \eqref{NullDecay02} was used in Klainerman~\cite{Kla86},
and the usage of $S$ in \eqref{NullDecay02} makes it difficult to treat the null form $Q_0$
included in coupled systems of the wave and Klein-Gordon equations,
and this is the reason why the notion of the strong null condition was introduced in \cite{Geo90}.

Before we proceed to our new estimate for $Q_0$,
we introduce another kind of known estimate for the null forms here. 
If we only use \eqref{frame02}, then we find
that the left-hand side of \eqref{NullRep01} is bounded from above by
$C\jb{r}^{-1}\left(\left|%\widetilde{Z} 
Z' \varphi\right|\, |\pa\psi|
{}+|\pa\varphi|\,\left|%\widetilde{Z} 
Z' \psi\right|
\right)$, where $Z':=(\Omega, \pa)$.
Hence, rewriting $\pa_+$ as
\begin{equation}
\label{frame05}
\pa_+=r^{-1}\left(S-(t-r)\pa_t \right),
\end{equation}
we get
\begin{equation}
\label{YokNull}
|Q_0(\varphi, \psi)|
\le
C \jb{r}^{-1}
\left(
\left|%\widetilde{\Gamma} 
\Gamma' \varphi\right|\,\left|\pa \psi\right|
{}+\left|\pa\varphi\right|\,\left|%\widetilde{\Gamma}
\Gamma' \psi\right|
{}+\jb{t-r}\left|\pa\varphi\right|\, \left|\pa\psi\right|
\right),
\end{equation}
where 
$\Gamma':=(S, \Omega, \pa)$.
Similar estimate can be obtained for $Q_{ab}$
in view of \eqref{NullRep04}.
These estimates are used in the study of
systems of wave equations with multiple speeds
because
$L$ is incompatible with such systems
(see Hoshiga-Kubo \cite{Hos-Kub00} and Yokoyama \cite{Yok00} for example). 
As we have mentioned in the previous section, 
the estimate \eqref{YokNull} is the point where
$S$ comes in the arguments of \cite{Yok00}, \cite{Kub-Yok01}
and \cite{Kat04:02}, though the weighted $L^\infty$--$L^\infty$ estimates (cf.~Lemma \ref{DecaySolWave}) are free of $S$.

% Let $w$ be a solution to $\dal w=\Phi$.
In Katayama-Kubo~\cite{Kat-Kub08:01}, 
for %$\pa_+ w$ 
the $\pa_+$-derivative of the solution to the wave equation,
%instead of rewriting $\pa_+$ by some vector fields,
a weighted $L^\infty$--$L^\infty$ estimate
with a better decay factor than \eqref{KubYokD} 
is directly obtained through an explicit expression of the solution
(without rewriting $\pa_+$ by the other vector fields),
%$w$ to the wave equation, 
and the null forms are treated using only 
$Z'=(\Omega, \pa)$ (see also \cite{Kat-Kub08:03}).
We can adopt this approach in \cite{Kat-Kub08:01}
to systems of the wave and Klein-Gordon equations
because the required vector fields $\Omega$ and $\pa$ are admissible.
However we take another approach here since we can use the vector field $L$;
motivated by \eqref{frame05}, we rewrite $\pa_+$ as
\begin{equation}
\label{frame06}
\pa_+=(t+r)^{-1} \left(2L_r+(t-r)\pa_t-(t-r)\pa_r\right).
\end{equation}
Then, by \eqref{NullRep01}, \eqref{NullRep03} and \eqref{frame06},
we obtain
\begin{equation}
\label{NullDecay03}
|Q_0(\varphi, \psi)|\le C\jb{t+r}^{-1} \left(|Z\varphi|\,|\pa \psi|+|\pa\varphi|\,|Z\psi|+\jb{t-r}|\pa\varphi|\,|\pa\psi|\right)
\end{equation}
at $(t,x)\in [0, \infty)\times \R^3$.

%%%%%
For any multi-index $\alpha$, we can easily check that 
$Z^\alpha Q_0(\varphi, \psi)$ can be written as a linear combination of 
the null forms $Q_0(Z^\beta \varphi, Z^\gamma \psi)$
and $Q_{cd}(Z^\beta \varphi, Z^\gamma \psi)$ with 
$|\beta|+|\gamma|\le |\alpha|$ and $0\le c<d\le 3$.
The same is true for $Z^\alpha Q_{ab}(\varphi, \psi)$ ($0\le a<b\le 3$).
Therefore, \eqref{NullDecay01} and \eqref{NullDecay03} yield the following:
%%%%%%%%%%%%%%%%%%%%%%%%%%%%%%%%%%
\begin{lemma}\label{SKnull}
Let $k$ be a nonnegative integer, 
and let $Q$ be one of the null forms $Q_0$ and $Q_{ab}$ with $0\le a<b\le 3$.
Then we have
\begin{align*}
\jb{t+|x|}|Q(\varphi, \psi)|_k
\le & C (|\varphi|_{[k/2]+1} |\pa \psi|_k+|\varphi|_{k+1} |\pa \psi|_{[k/2]})\\
& {}+C(|\pa \varphi|_{[k/2]} |\psi|_{k+1}+|\pa \varphi|_k |\psi|_{[k/2]+1})\nonumber\\
& +C\jb{t-|x|}(|\pa \varphi|_{[k/2]} |\pa \psi|_k
{}+|\pa \varphi|_k |\pa \psi|_{[k/2]})
\nonumber
\end{align*}
at $(t,x)\in [0,\infty)\times \R^3$
for any smooth functions $\varphi$ and $\psi$.
Here $C$ is a positive constant depending only on $k$, 
$|\cdot|_s$ is given by \eqref{InvariantNorm}
for a nonnegative integer $s$, and $[m]$ 
denotes the largest integer not exceeding the number $m$. 
\end{lemma}

%%%%%%%%%%%%%%%%%%%%%%%%%%%%%%%%%%%%%%%%%%%
\section{Proof of Theorem \ref{Global}}
%%%%%%%%%%%%%%%%%%%%%%%%%%%%%%%%%%%%%%%%%%%
In this section, we will prove Theorem \ref{Global}.

Suppose that all the assumptions in Theorem \ref{Global} are fulfilled.
The classical theory for nonlinear hyperbolic equations
implies %We can easily obtain 
the local existence of the classical solutions
to \eqref{OurSystem}--\eqref{Data} for small $\ve$.
%under the conditions \eqref{quasi-linear} and \eqref{symmetricity}.
Moreover, we see that the solution $u$ exists as long as
$\sum_{|\alpha|\le 2} \|\pa^\alpha u(t,\cdot)\|_{L^\infty(\R^3)}$ stays finite
(see H\"ormander \cite{Hoe97} for instance).
Hence what we need 
for the proof of Theorem \ref{Global}
is such an {\it a priori} estimate to guarantee the boundedness of $\sum_{|\alpha|\le 2}\|\pa^\alpha u(t,\cdot)\|_{L^\infty(\R^3)}$.
Let $u=(u_i)_{1\le i\le N}=(v, w)$ be the local solution
to \eqref{OurSystem}--\eqref{Data} for $0\le t<T_0$
with some $T_0>0$, where $v$ and $w$
are given by \eqref{division}.
If both ${\mathcal I}_1$ and ${\mathcal I}_2$ in the condition (b) are non-empty,
without loss of generality we may assume
that ${\mathcal I}_1=\{1, \ldots, N_3\}$, and ${\mathcal I}_2=\{N_3+1, \ldots, N_2\}$
with some positive integer $N_3$.
Correspondingly, we write
$$
w=(w_k)_{1\le k\le N_2}=\bigl( (w_k^{\rm (i)})_{1\le k\le N_3}, (w_k^{\rm (ii)})_{1\le k\le N_4} \bigr)
=\bigl(w^{\rm (i)}, w^{\rm (ii)} \bigr),
$$
where $N_4=N_2-N_3$. If ${\mathcal I}_2$ (resp.~${\mathcal I}_1$) is empty, 
then we put $w^{\rm (i)}=w$ (resp.~$w^{\rm (ii)}=w$), and 
$w^{\rm (ii)}$ (resp.~$w^{\rm (i)}$) should be neglected in what follows. 

For a nonnegative integer $\sigma$, and a positive constant $p$, we define
\begin{align*}
d_{\sigma, p}(t,x)=& 
\jb{t+|x|}^{3/2} |v(t,x)|_{\sigma+2}+\jb{x}\jb{t-|x|} |\pa w(t,x)|_{\sigma+1} \\
& {}+\jb{t+|x|}\left({\mathcal W_0}(t,|x|) |w^{\rm (i)}(t,x)|_{\sigma+2}
{}+W_-(t,|x|)^{1-p}|w^{\rm (ii)}(t,x)|_{\sigma+2}\right)
\end{align*}
for $(t,x)\in [0,T_0)\times \R^3$, where $|\cdot|_s$, ${\mathcal W}_0$,
and $W_-$ are given by \eqref{InvariantNorm}, \eqref{KY-weight}, and \eqref{WeakWei},
respectively.

For a smooth function $\varphi=\varphi(x)$ and a nonnegative integer $s$, we set
\begin{align*}
\|\varphi\|_{X^s}=& \sum_{|\alpha|\le s}\left(\left(\int_{\R^3} 
|\jb{x}^{s+2}\pa_x^\alpha \varphi(x)|^2dx\right)^{1/2}
+\left(\int_{\R^3} |\jb{x}^{s}\pa_x^\alpha\varphi(x)|^{6/5}dx\right)^{5/6}
\right).
\end{align*}
Note that the Sobolev embedding theorem implies that
\begin{align*}
\sup_{x\in \R^3} \sum_{|\alpha|\le s-2}
\jb{x}^{|\alpha|+4} |\pa_x^\alpha \varphi(x)|\le \sup_{x\in \R^3} \sum_{|\alpha|\le s-2}
\jb{x}^{s+2} |\pa_x^\alpha \varphi(x)|
\le C_s\|\varphi\|_{X^s}
\end{align*}
for $s\ge 2$, where $C_s$ is a positive constant depending only on $s$.

Our aim here is to show the following:
\begin{proposition}\label{bootstrap}
Fix some $\sigma\ge 19$, and $0<p<1/100$, say.
Suppose that all the assumptions in Theorem \ref{Global} are fulfilled.
Assume that $\|f\|_{X^{2\sigma+1}}+\|g\|_{X^{2\sigma}}\le M_0$ with some positive constant $M_0$.
Let $u=(v,w)$ be the local solution to \eqref{OurSystem}--\eqref{Data} 
for $0\le t<T_0$.
%, and suppose $0<T\le T_0$.
Then there exists a positive constant $A_0=A_0(M_0)$ having the following property:
For any $A\ge A_0$, there exists a positive constant $\ve_0=\ve_0(A)$ such that
\begin{equation}
\label{boot01}
 \sup_{0\le t<T} \|d_{\sigma,p}(t, \cdot)\|_{L^\infty(\R^3)} \le A\ve
\end{equation}
implies 
\begin{equation}
\label{boot-finale}
 \sup_{0\le t<T} \|d_{\sigma,p}(t,\cdot)\|_{L^\infty(\R^3)}\le \frac{A}{2}\ve,
\end{equation}
provided that $0<\ve\le \ve_0$ and $0<T\le T_0$.
Here $A_0$ and $\ve_0$ are independent of %$T$ and 
$T_0$.
\end{proposition}
Once Proposition \ref{bootstrap} is established, by the
continuity argument (or the bootstrap argument),
we find that $\|d_{\sigma,p}(t,\cdot)\|_{L^\infty(\R^3)}$ stays
bounded as long as the solution exists, provided that $\ve$ is small enough.
Indeed, suppose that $f$ and $g$ belong to $C^\infty_0(\R^3;\R^N)$ at first.
Then, taking the support of $u$ into account, we
see that $\|d_{\sigma,p}(t,\cdot)\|_{L^\infty(\R^3)}$ is continuous in $t$.
We choose a large constant $A(\ge A_0)$
to satisfy $\|d_{\sigma,p}(0,\cdot)\|_{L^\infty(\R^3)}< A\ve$.
Then we see that \eqref{boot01} is true for some small $T$.
Let $T_*(>0)$ be the supremum of $T\left(\in (0, T_0)\right)$ for which \eqref{boot01} holds.
If $\ve\in (0, \ve_0]$, then by \eqref{boot-finale} and the continuity of $\|d_{\sigma, p}(t,\cdot)\|_{L^\infty}$,
we conclude that $T_*=T_0$ (otherwise we meet a contradiction). In other words,
if $u$ is the local solution for $0\le t<T_0$,
then we have $\sup_{0\le t<T_0} \|d_{\sigma, p}(t,\cdot)\|_{L^\infty}\le A\ve$,
provided that $\ve\in (0, \ve_0]$.
We see that the same is true for general $f,g\in{\mathcal S}(\R^3;\R^N)$ through the
approximation by $C^\infty_0$-functions.
This {\it a priori} estimate implies Theorem \ref{Global} immediately.

Now we are going to prove Proposition \ref{bootstrap}.
We assume that \eqref{boot01} holds.
In the following, various positive constants, being independent of $A(>0)$, $\ve(\le 1)$,
$T(>0)$, and $M_0$, are indicated just by the same letter $C$. Thus the practical value of $C$ may change line by line. 
Similarly $C_*$ stands for various positive constants depending only on $M_0$ 
and the bounds for finite numbers of derivatives of $F$ in a small neighborhood
of the origin $(\xi,\xi',\xi'')=(0,0,0)$.
We always assume that $\ve$ is small enough to
satisfy $A\ve\le 1$, say.

First we remark that for any nonnegative integer $s$, there exists a positive constant $C_s$ such that
\begin{equation}
\label{CommZD2}
C_s^{-1} |\pa \varphi(t,x)|_s\le \sum_{|\alpha|\le s} |\pa Z^\alpha \varphi(t,x)|
\le C_s |\pa\varphi(t,x)|_s
\end{equation}
holds for any smooth function $\varphi$, because of \eqref{CommZD}.
We also note that
we have
 \begin{equation}
 \label{KY}
  \jb{x}^{-1}\jb{t-|x|}^{-1}\le C \jb{t+|x|}^{-1}{W}_-(t,|x|)^{-1},\quad (t,x)\in [0,\infty)\times \R^3.
 \end{equation}
We fix some small and positive constant $\delta$. 
Then we have
\begin{equation}
\label{boot01c}
 {\mathcal W}_{0}(t,|x|)^{-1}\le C\jb{t+|x|}^\delta \jb{t-|x|}^{-\delta},\quad (t,x)\in [0,\infty)\times \R^3.
\end{equation}
We will use \eqref{KY} and \eqref{boot01c} repeatedly in the following. 
Note that \eqref{boot01}, \eqref{KY}, and \eqref{boot01c} yield
\begin{equation}
|u(t,x)|_{\sigma+2}\le CA\ve \jb{t+|x|}^{-1+\delta}\jb{t-|x|}^{-\delta},\quad (t,x)\in [0,T)\times \R^3.
\label{boot01d}
\end{equation}

The proof of Proposition \ref{bootstrap} is divided into several steps.
%%%%%%%%%%%%%%%%
\smallskip

{\it Step 1: Energy Estimate.}
Let $0<\lambda<p/4$. In this step, we are going to prove that
\begin{equation}
\label{Ene01}
\sup_{0\le t<T} (1+t)^{-\lambda} \left(\|v(t)\|_{2\sigma}+\|w^{\rm (ii)}(t)\|_{2\sigma}+\|\pa u(t)\|_{2\sigma}\right)\le C_*\ve
\end{equation}
holds for small $\ve$, where $\|\cdot\|_s$ is given by \eqref{InvariantNorm}.
The difficulty here is the lack of a natural estimate for $\|w^{\rm (i)}(t)\|_{2\sigma}$ (cf.~Lemma~\ref{EnergyHyp}).
To overcome this difficulty, we will use the following lemma
that is easily obtained from the definition of $Z$ and \eqref{CommZD2}:
\begin{lemma}
For any $s\ge 1$, there exists a positive constant $C=C(s)$ such that we have
\begin{equation}
\label{KKTech}
|\varphi(t,x)|_s\le C\left( |\varphi(t,x)|+\jb{t+|x|}|\pa \varphi(t,x)|_{s-1} \right)
\end{equation}
for any smooth function $\varphi=\varphi(t,x)$.
\end{lemma}
In fact, for $s\ge 1$, we have
\begin{align*}
|\varphi(t,x)|_s \le & C\Bigl(|\varphi(t,x)|+\sum_{|\alpha|=1}\sum_{|\beta|\le s-1} |Z^\alpha Z^\beta \varphi(t,x)|\Bigr)\\
\le & C\Bigl(|\varphi(t,x)|+\jb{t+|x|}\sum_{|\beta|\le s-1} |\pa (Z^\beta \varphi)(t,x)|\Bigr),
\end{align*}
which leads to \eqref{KKTech}, thanks to \eqref{CommZD2}.

Now we start the proof of \eqref{Ene01}.
Let $|\alpha|=s\le 2\sigma$. We set
\begin{equation}
\label{Right01}
{F}_{i, \alpha}
=Z^\alpha \left\{F_i(u, \pa u, \pa_x\pa u)\right\}-\sum_{j,k,a} \gamma_{ka}^{\, ij}(u,
\pa u)\pa_k\pa_a(Z^\alpha u_j),
\end{equation}
where $\gamma=(\gamma_{ka}^{\, ij})$ is from \eqref{quasi-linear}.
Then we have
\begin{equation}
\label{HigherSys01}
\left(\dal+m_i^2\right)(Z^\alpha u_i)-\sum_{j,k,a}\gamma_{ka}^{\, ij}(u, \pa u)
\pa_k\pa_a (Z^\alpha u_j)
={F}_{i,\alpha}, \ 1\le i\le N.
\end{equation}
Note that we have $\left|[Z^\alpha, \pa_k\pa_a]u_j\right|\le C |\pa u|_s$ by
\eqref{CommZD}. Hence, in view of \eqref{quasi-linear}, \eqref{CommZD}, and \eqref{Right01}, 
we obtain from the condition (b--i), \eqref{boot01}, \eqref{KY}, and \eqref{boot01d}
that
\begin{align}
\left| 
{F}_{i, \alpha}\right|\le & 
C \left(|v|_{[s/2]}+|w^{\rm (ii)}|_{[s/2]}+|\pa u|_{[s/2]+1}\right)(|v|_{s}+|w^{\rm (ii)}|_s+|\pa u|_{s})
\nonumber\\
& 
{}+C|u|_{[s/2]+2}^2
 (|u|_{s}+|\pa u|_{s}) \nonumber
 \\
 \le & CA\ve \jb{t+|x|}^{-1}(|v|_{s}+|w^{\rm (ii)}|_s+|\pa u|_{s})
\nonumber\\
& 
{}+CA^2\ve^2 \jb{t+|x|}^{-2+2\delta}\jb{t-|x|}^{-2\delta} \nonumber\\
 & \qquad \times
     \left(A\ve\jb{t+|x|}^{-1+\delta}\jb{t-|x|}^{-\delta}+\jb{t+x}|\pa u|_{s-1}+|\pa u|_{s}\right) 
\label{Right02}
\end{align}
at $(t,x)\in [0, T)\times \R^3$. Here we have also used \eqref{KKTech} to estimate $|u|_s$ for $s\ge 1$.
Thus the term $\jb{t+|x|}|\pa u|_{s-1}$ on the right-hand side of \eqref{Right02} should be neglected when $s=0$.
Since we have
$$
\left\|\jb{t+|\cdot|}^{-3+3\delta}\jb{t-|\cdot|}^{-3\delta}\right\|_{L^2(\R^3)}\le C (1+t)^{-3/2}
$$
for $\delta<1/6$, \eqref{Right02} yields
\begin{align}
\nonumber
\left\| 
{F}_{i, \alpha}
\right\|_{L^2(\R^3)}\le & CA\ve (1+t)^{-1}(\|v\|_{s}+\|w^{\rm (ii)}\|_{s}+\|\pa u\|_{s})\\
& {}+CA^2\ve^2 (1+t)^{-1+2\delta}\|\pa u\|_{s-1}
+CA^3\ve^3(1+t)^{-3/2},
\label{Right03}
\end{align}
where the term $CA^2\ve^2(1+t)^{-1+2\delta}\|\pa u\|_{s-1}$ should be neglected when $s=0$.
In view of the condition (b--i), we also obtain from \eqref{boot01}, \eqref{KY}, and \eqref{boot01d} that
\begin{equation}
\label{coeff01}
|\gamma|_1\le C(|v|_1+|w^{\rm (ii)}|_1+|\pa u|_1+|u|_1^2+|\pa u|_1^2)\le CA\ve \jb{t+|x|}^{-1}
\end{equation}
at $(t,x)\in [0,T)\times \R^3$.
Because of \eqref{symmetricity} and \eqref{coeff01}, we can apply Lemma \ref{EnergyHyp}
to \eqref{HigherSys01} for small $\ve$.

For $N_3+1\le k\le N_2$ and $0\le a\le 3$, let ${\mathcal G}_{k,a}$ be from the condition (b--ii).
Because of \eqref{CommZD}, we get
\begin{equation}
\label{DivHigher01}
\dal \left(Z^\alpha w_j^{\rm (ii)} \right)=\sum_{a,b=0}^3 \sum_{|\beta|\le |\alpha|}
C^{\alpha\beta}_{ab} \pa_b \left(Z^\beta {\mathcal G}_{N_3+j, a}(u,\pa u)\right),\quad 1\le j\le N_4
\end{equation}
with appropriate constants $C^{\alpha\beta}_{ab}$.
Remember that each ${\mathcal G}_{k,a}^{\rm (q)}$ is independent of $w^{\rm (i)}$ itself.
Thus, going similar lines to \eqref{Right02} and \eqref{Right03}, we get
\begin{align}
\nonumber
\left\|
{\mathcal G}_{k,a}(u,\pa u)
\right\|_{s}\le & CA\ve (1+t)^{-1}(\|v\|_{s}+\|w^{\rm (ii)}\|_{s}+\|\pa u\|_{s})\\
& {}+CA^2\ve^2 (1+t)^{-1+2\delta}\|\pa u\|_{s-1}
+CA^3\ve^3(1+t)^{-3/2}
\label{Right03'}
\end{align}
for $s\le 2\sigma$. As before, the term including
$\|\pa u\|_{s-1}$ on the right-hand side should be neglected when $s=0$.

We put
$$
E_s(t)=\|v(t)\|_{s}+\|w^{\rm (ii)}(t)\|_s+\|\pa u(t)\|_{s}
$$
for $s\ge 0$.
Applying Lemma \ref{EnergyHyp}
to \eqref{HigherSys01} with $|\alpha|=s=0$,
and applying Lemma \ref{DivEne} to \eqref{DivHigher01} with $|\alpha|=s=0$,
we obtain from \eqref{Right03}, \eqref{coeff01}, and \eqref{Right03'}
that
$$
E_0(t)
\le C_*\ve+CA^3\ve^3
{}+CA\ve \int_0^t (1+\tau)^{-1}
E_0(\tau) d\tau.
$$
The Gronwall lemma yields
\begin{equation}
E_0(t)
\le \left(C_*\ve+CA^3\ve^3\right)(1+t)^{CA\ve} \le C_*\ve (1+t)^{CA\ve},
\label{STAT}
\end{equation}
provided that $\ve$ is small enough to satisfy $A^3\ve^2\le 1$.
Starting with \eqref{STAT}, we can inductively obtain
\begin{equation}
\label{Ene02}
E_s(t)
\le C_{*,s} \ve (1+t)^{2s\delta+CA\ve}
\end{equation}
for $0\le s\le 2\sigma$, where $C_{*,s}$'s are positive constants depending on
$s$, $M_0$ and the nonlinearity $F$. 
In fact, if \eqref{Ene02} with $s$ replaced by $s-1$ is true for some $s\ge 1$, then 
applying Lemmas~\ref{EnergyHyp} and \ref{DivEne} to \eqref{HigherSys01} and
\eqref{DivHigher01} with $|\alpha|=s$, respectively,
and using \eqref{Right03}, \eqref{coeff01} and
\eqref{Right03'}, we obtain
\begin{align*}
E_s(t)
\le & C_*\ve+C\left(A^3\ve^3+(2s\delta+CA\ve)^{-1}C_{*, s-1}A^2\ve^3(1+t)^{2s\delta+CA\ve}\right)
\\
& {}+CA\ve \int_0^t (1+\tau)^{-1}
E_s(\tau) d\tau, \nonumber
\end{align*}
and the Gronwall lemma leads to \eqref{Ene02}.

Finally, 
we obtain \eqref{Ene01} from \eqref{Ene02} with $s=2\sigma$,
provided that $\delta$ in \eqref{boot01c} is chosen to satisfy $4\sigma\delta\le \lambda/2$, and 
$\ve$ is small enough to satisfy %$A^3\ve^2\le 1$ and 
$CA\ve\le \lambda/2$.
%%%%%%%%%%%%%%%%%%%%%%%%%%%%%%%%%%%%%
\smallskip

{\it Step 2: Decay Estimates, Part 1.}
By Lemma \ref{KlaSobolev} and \eqref{Ene01}, we get
\begin{equation}
\label{Decay01}
\jb{x}\left(|v(t,x)|_{2\sigma-2}+|w^{\rm (ii)}(t,x)|_{2\sigma-2}+|\pa u(t,x)|_{2\sigma-2}\right)\le C_*\ve (1+t)^\lambda
\end{equation}
for $(t,x)\in [0,T)\times \R^3$.

Similarly to \eqref{Right02}, we get
\begin{align}
\nonumber
|F_i|_{s} \le & 
C \bigl(|v|_{\sigma}+|w^{\rm (ii)}|_{\sigma}+|\pa u|_{\sigma+1}\bigr) 
\bigl(|v|_{s}+|w^{\rm (ii)}|_{s}+|\pa u|_{s+1}\bigr)
\\& 
{}+C|u|_{\sigma+2}^2 \left(|u|_{s}+|\pa u|_{s+1}\right) \nonumber
\\
\le & CA\ve\left(
                  \jb{t+|x|}^{-3/2}+\jb{t+|x|}^{-1}{W}_-(t,|x|)^{-1+p}
                 \right)
%\nonumber\\
%& \qquad\qquad \times
 \bigl( |v|_{s}+|w^{\rm (ii)}|_{s}+|\pa u|_{s+1}\bigr) \nonumber\\
& {}+CA^2\ve^2 \jb{t+|x|}^{-2+2\delta}\jb{t-|x|}^{-2\delta} 
 \bigl(|v|_{s}+|w^{\rm (ii)}|_{s}+|\pa u|_{s+1}+|w^{\rm (i)}|_s\bigr)
       \label{Right04}                    
\end{align}
for $s\le 2\sigma$.
For $\rho\ge 0$ and a nonnegative integer $s$, we set
\begin{equation}\label{DefM}
 M^{\rm (i)}_{\rho, s}=\sup_{(t,x)\in [0,T)\times \R^3} \jb{t+|x|}^{1-\rho}{\mathcal W}_0(t,|x|)|w^{\rm (i)}(t,x)|_s.
\end{equation}
Then, using \eqref{Decay01} and \eqref{Right04}, we get
\begin{align}
\nonumber
\jb{x}|F_i|_{2\sigma-3}
\le & C_*A\ve^2 \jb{t+|x|}^{-1+\lambda}{W}_-(t,|x|)^{-1/2}
\\
& {}+C_*A^2\ve^3 \jb{t+|x|}^{-2+2\delta+\lambda}\jb{t-|x|}^{-2\delta}
\nonumber\\
& {}+CA^2\ve^2M^{\rm (i)}_{\lambda+(1/2), 2\sigma-3} \jb{t+|x|}^{-(3/2)+3\delta+\lambda}\jb{t-|x|}^{-3\delta}
\nonumber\\
\le & \left(C_*A\ve^2+CA^2\ve^2M^{\rm (i)}_{\lambda+(1/2), 2\sigma-3}\right)\jb{t+|x|}^{-(1/2)+\lambda-\mu}{W}_-(t,|x|)^{-1+\mu}, \nonumber
\end{align}
where $\mu$ is a small and positive constant.
Hence, by Lemma~\ref{Asakura}, and also by \eqref{KubYok} of Lemma~\ref{DecaySolWave} with $(\rho,\kappa)=(\lambda+(1/2), 1)$, we get
$$
M^{\rm (i)}_{\lambda+(1/2), 2\sigma-3}\le 
C_*\left(\ve+A\ve^2\right)+CA^2\ve^2M^{\rm (i)}_{\lambda+(1/2),2\sigma-3}.
$$
Therefore, if $\ve$ is small enough to satisfy 
$CA^2\ve^2\le 1/2$, then we obtain
\begin{align}
\sup_{(t,x)\in[0,T)\times \R^3}\jb{t+|x|}^{-\lambda+(1/2)}{\mathcal W}_0(t,|x|)|w^{\rm (i)}(t,x)|_{2\sigma-3}
=& M^{\rm (i)}_{\lambda+(1/2), 2\sigma-3} 
 \nonumber\\
\le & C_*\ve. 
\label{Decay02}
\end{align}

Using \eqref{Decay02}, we have
\begin{equation}
 \jb{x}|F_i|_{2\sigma-3}\le C_*(A\ve^2+A^2\ve^3)\jb{t+|x|}^{-(1/2)+\lambda-\mu}{W}_-(t,|x|)^{-1+\mu}.
\label{KY02}
\end{equation}
Hence, similarly to \eqref{Decay02}, Lemma \ref{Asakura} and \eqref{KubYokD} of
Lemma~\ref{DecaySolWave} yield
\begin{equation}
\label{Decay03}
\sup_{(t,x)\in [0,T)\times \R^3} \jb{t+|x|}^{-\lambda-(1/2)}\jb{x}\jb{t-|x|}|\pa w(t,x)|_{2\sigma-4}
\le C_*\ve,
\end{equation}
provided that $\ve$ is small enough.
Going similar lines to \eqref{Right04}--\eqref{KY02}, we get
$$
 \jb{x} |{\mathcal G}_{k,a}(u,\pa u)|_{2\sigma-3}\le C_*(A\ve^2+A^2\ve^3)\jb{t+|x|}^{-(1/2)+\lambda-\mu}{W}_-(t,|x|)^{-1+\mu},
$$
and applying Lemma~\ref{DivDecay} to \eqref{DivHigher01} with $|\alpha|\le 2\sigma-4$,
we get
\begin{equation}
\label{Decay03'}
\sup_{(t,x)\in [0,T)\times \R^3} \jb{t+|x|}^{-\lambda-(1/2)}\jb{x}\jb{t-|x|}|w^{\rm (ii)}(t,x)|_{2\sigma-4}
\le C_*\ve,
\end{equation}
provided that $\ve$ is small enough.

Using  \eqref{Ene01}, 
we obtain from \eqref{Right01}, \eqref{Right02}, and \eqref{coeff01} that 
\begin{align}
\nonumber
\|\jb{t+|\cdot|}|F_i|_{2\sigma-1}\|_{L^2}\le & C_*A\ve^2 (1+t)^\lambda
{}+C_*A^2\ve^3(1+t)^{\lambda+2\delta}\\
& {}+C A^3\ve^3
\left\|
\jb{t+|\cdot|}^{-2+3\delta}\jb{t-|\cdot|}^{-3\delta}
\right\|_{L^2}
\le
C_*A\ve^2 (1+t)^{2\lambda}
\nonumber
\end{align}
for sufficiently small $\delta$. %, and $0<\ve<A^{-1}$.
Hence Lemma \ref{DecayKleinGordon} leads to
\begin{align}
\label{Decay04a}
\jb{t+|x|}^{3/2}|v(t,x)|_{2\sigma-5}\le & C_*\left(\ve+A\ve^2 \sum_{j=0}^\infty \sup_{\tau\in (0,t)}
\chi_j(\tau)(1+\tau)^{2\lambda}\right).
\end{align}
Let $2^{J-1}\le t< 2^{J}$ with some nonnegative integer $J$. Then we have
\begin{align*}
\sum_{j=0}^\infty \sup_{\tau\in (0,t)}
\chi_j(\tau)(1+\tau)^{2\lambda}=& \sum_{j=0}^{J} \sup_{\tau\in (0,t)}
\chi_j(\tau)(1+\tau)^{2\lambda}\\
\le & \sum_{j=0}^{J} 2^{2(j+2)\lambda}= \frac{2^{4\lambda}
\left(2^{2\lambda(J+1)}-1\right)}{2^{2\lambda}-1}
\le C(1+t)^{2\lambda}. \nonumber
\end{align*}
A similar estimate for $0\le t<1$ is trivially obtained.
Now \eqref{Decay04a} leads to
\begin{equation}
\label{Decay04}
\jb{t+|x|}^{(3/2)-2\lambda}|v(t,x)|_{2\sigma-5}\le C_*\ve,
\end{equation}
provided that $\ve$ is small enough.
%%%%%%%%%%%%%%%%%%%%%%%%%%%%%%%%%%%%%%%%%%%

{\it Step 3: Decay Estimates, Part 2.}
We make use of the detailed structure of the nonlinearity
from now on.
%%%%%%%%%%%%%%%%%%%%%%%%%%%%%%%%%%%%%%%%%%%%%%%%%%
Recall that $\Phi^{\rm (W)}$ is given by \eqref{FIW}
for a smooth function $\Phi=\Phi(\xi, \xi',\xi'')$.
We also define
\begin{align}
\label{FIK}
\Phi^{\rm (K)}(\eta,\eta',\eta''):=& \Phi^{\rm (q)}\bigl((\eta,\zeta), (\eta',\zeta'), (\eta'',\zeta'')\bigr)\Bigr|_{(\zeta,\zeta',\zeta'')=(0,0,0)},\\
\label{FIKW}
\Phi^{\rm (KW)}(\xi, \xi',\xi''):=& 
\Phi^{\rm (q)}\bigl((\eta,\zeta), (\eta',\zeta'), (\eta'',\zeta'')\bigr)\\
  & {}-\Phi^{\rm (K)}(\eta,\eta', \eta'')-\Phi^{\rm (W)}(\zeta,\zeta',\zeta''), \nonumber\\
\label{FIH}
\Phi^{\rm (H)}(\xi,\xi',\xi''):=& \Phi(\xi,\xi', \xi'')-\Phi^{\rm (q)}(\xi, \xi', \xi''),
\end{align}
where $(\xi, \xi', \xi'')=\bigl((\eta,\zeta), (\eta',\zeta'), (\eta'',\zeta'')\bigr)$
as before 
(the letters ``K'', ``W'', ``H'' in this
 notation stand for ``Klein-Gordon'', ``wave'',
``higher nonlinearity'', respectively).

%%%%%%%%%%%%%%%%%%%%%%%%%%%%%%%%%%%%%%%%%%%%%%%%%%%%%%
Since $\bigl(F_i^{\rm (W)}\bigr)_{N_1+1\le i\le N}$
satisfies the null condition, 
by \eqref{NullFormNullCond}, Lemma \ref{SKnull}, and \eqref{boot01}, we get
\begin{align}
\nonumber
\bigl|F_i^{\rm (W)}\bigr|_s \le & C\jb{t+|x|}^{-1} \left(|w|_{[s/2]+2}|\pa w|_{s+1}+|\pa w|_{[s/2]+1}|w|_{s+1}
\right) \\
&{}+C\jb{t+|x|}^{-1}\jb{t-|x|}|\pa w|_{[s/2]+1}|\pa w|_{s+1}\nonumber\\
\le &CA\ve\jb{t+|x|}^{-2+\delta}
\jb{t-|x|}^{-\delta}
|\pa w|_{s+1}\nonumber\\
&{}+CA\ve \jb{x}^{-1}\jb{t+|x|}^{-1}\left(\jb{t-|x|}^{-1}|w|_{s+1}+|\pa w|_{s+1}\right)
\label{Right11}
\end{align}
at $(t,x)\in [0,T)\times \R^3$
for $N_1+1\le i\le N$, provided that $s\le 2\sigma$.

For $1\le i\le N$, it is easy to see
that each $F_i^{\rm (K)}(v, \pa v, \pa_x\pa v)$ is a linear combination
of $(\pa^\alpha v_j)(\pa^\beta v_k)$ with $|\alpha|, |\beta|\le 2$, and
$1\le j,k \le N_1$. %On the other hand, each 
Similarly we can see that each $F_i^{\rm (KW)}(u, \pa u, \pa_x\pa u)$ is a linear combination of
$(\pa^\alpha v_j)(\pa^\beta w_k)$
and $(\pa^\alpha v_j)w^{\rm (ii)}_l$ with $|\alpha| \le 2$,
$1\le |\beta|\le 2$,  
$1\le j\le N_1$, $1\le k\le N_2$, and $1\le l\le N_4$, because of the condition (b--i).
Therefore \eqref{boot01} yields
\begin{align}
\label{Right12}
\bigl|F_i^{\rm (K)}\bigr|_{s}
        \le & C|v|_{[s/2]+2} |v|_{s+2}\le CA\ve\jb{t+|x|}^{-3/2}|v|_{s+2},\\
\bigl|F_i^{\rm (KW)}\bigr|_s
        \le & C\left(|w^{\rm (ii)}|_{[s/2]}+|\pa w|_{[s/2]+1}\right)|v|_{s+2}
{}+|v|_{[s/2]+2}
\left(|w^{\rm (ii)}|_s+|\pa w|_{s+1}\right)
\nonumber\\
\le & CA\ve \jb{t+|x|}^{-1}W_-(t, |x|)^{-1+p} |v|_{s+2}\nonumber\\
&{}+CA\ve\jb{t+|x|}^{-3/2}\left(|w^{\rm (ii)}|_s+|\pa w|_{s+1}\right)
\label{Right13}
\end{align}
at $(t,x)\in [0, T)\times \R^3$ for $1\le i\le N$, provided that $s\le 2\sigma$.

Since we have 
$F_i^{\rm (H)}(u, \pa u, \pa_x\pa u)=O\left(|u|^3+|\pa u|^3+|\pa_x\pa u|^3\right)$,
we get
\begin{align}
\nonumber
\bigl|F_i^{\rm (H)}|_s\le & C|u|_{[s/2]+2}^2\left(|v|_{s+2}+|w|_s+|\pa w|_{s+1}\right)
\\
\le & CA^2\ve^2 \jb{t+|x|}^{-2+2\delta}\jb{t-|x|}^{-2\delta}
\nonumber\\
& \qquad\qquad \times \bigl(|v|_{s+2}+|w^{\rm (ii)}|_s+|\pa w|_{s+1}{}+|w^{\rm (i)}|_s\bigr)
\label{Right14}
\end{align}
at $(t,x)\in [0, T)\times \R^3$ for $1\le i\le N$, provided that $s\le 2\sigma$.

Let $N_1+1\le i\le N$ in the following.
Using \eqref{Decay02}, \eqref{Decay03}, \eqref{Decay03'}, and \eqref{Decay04}, 
we obtain from \eqref{Right11}, \eqref{Right12}, and \eqref{Right13}
that
\begin{align}
\label{Right21}
\jb{x}\bigl| F_i^{\rm (W)} \bigr|_{2\sigma-7}
      \le & C_*A\ve^2\jb{t+|x|}^{-(3/2)+\delta+\lambda}
{W}_-(t,|x|)^{-1-\delta},
\\
\label{Right22}
\jb{x}\bigl|F_i^{\rm (K)}\bigr|_{2\sigma-7}\le & C_*A\ve^2 \jb{t+|x|}^{-2+2\lambda},\\
\label{Right23}
\jb{x}\bigl|F_i^{\rm (KW)}\bigr|_{2\sigma-7}\le & C_*A\ve^2 \jb{t+|x|}^{-1+\lambda}
{W}_-(t,|x|)^{-1}
\end{align}
at $(t,x)\in [0, T)\times \R^3$.
Similarly, by \eqref{Right14}, we get
\begin{align}
\nonumber
\jb{x}\bigl|F_i^{\rm (H)}\bigr|_{2\sigma-7}\le & C_*A^2\ve^3\jb{t+|x|}^{-(3/2)+2\delta+2\lambda}\jb{t-|x|}^{-1-2\delta}\\
                           &{}+CA^2\ve^2M^{\rm (i)}_{3\lambda, 2\sigma-7}\jb{t+|x|}^{-2+3\delta+3\lambda}
                               \jb{t-|x|}^{-3\delta}
                           \label{Right24}
\end{align}
at $(t,x)\in [0, T)\times \R^3$, where $M^{\rm (i)}_{3\lambda, 2\sigma-7}$ is given by
\eqref{DefM}.
Since we can choose $\delta$ as small as we wish,
\eqref{Right21} -- \eqref{Right24} lead to
\begin{equation}
\label{lalala}
\jb{x}|F_i|_{2\sigma-7}\le \left(C_*A\ve^2+CA^2\ve^2M^{\rm (i)}_{3\lambda, 2\sigma-7}\right)\jb{t+|x|}^{-1+3\lambda-\mu}
{W}_-(t,|x|)^{-1+\mu}
\end{equation}
for small $\mu>0$.
Therefore, by \eqref{KubYok} with $(\rho,\kappa)=(3\lambda,1)$, and by Lemma \ref{Asakura}, 
we obtain
$$
M^{\rm (i)}_{3\lambda, 2\sigma-7}\le C_*\left(\ve+A\ve^2\right)+CA^2\ve^2M^{\rm (i)}_{3\lambda, 2\sigma-7},
$$
which leads to
\begin{equation}
\label{Decay11}
\sup_{(t,x)\in [0, T)\times \R^3}\jb{t+|x|}^{1-3\lambda}{\mathcal W}_0(t,|x|)
|w^{\rm (i)}(t,x)|_{2\sigma-7} =M^{\rm (i)}_{3\lambda, 2\sigma-7} \le C_*\ve,
\end{equation}
provided that $\ve$ is sufficiently small.
Now \eqref{lalala} and \eqref{Decay11} yield
\begin{equation}
\label{KY10}
\jb{x}|F_i|_{2\sigma-7}\le C_*A\ve^2\jb{t+|x|}^{-1+3\lambda-\mu}
W_-(t,|x|)^{-1+\mu}.
\end{equation}
Hence Lemma \ref{Asakura} and \eqref{KubYokD} imply
\begin{equation}
\label{Decay12}
\sup_{(t,x)\in [0, T)\times \R^3}\jb{t+|x|}^{-3\lambda}\jb{x}\jb{t-|x|}
|\pa w(t,x)|_{2\sigma-8} \le C_*\ve,
\end{equation}
provided that $\ve$ is sufficiently small.

Now we are going to estimate $w^{\rm (ii)}$.
Suppose that $N_3+1\le k\le N_2$, and $0\le a\le 3$. 
Let ${\mathcal G}_{k,a}^{\rm (W)}(\zeta, \zeta')$, ${\mathcal G}_{k,a}^{\rm (K)}(\eta, \eta')$,
${\mathcal G}_{k,a}^{\rm (KW)}(\xi, \xi')$, and ${\mathcal G}_{k,a}^{\rm (H)}(\xi, \xi')$ be given by
\eqref{FIW}, \eqref{FIK}, \eqref{FIKW}, and \eqref{FIH}, respectively,
with $\Phi={\mathcal G}_{k,a}(\xi,\xi')$.
We can easily verify that ${\mathcal G}_{k,a}^{\rm (K)}$, ${\mathcal G}_{k,a}^{\rm (KW)}$, and ${\mathcal G}_{k,a}^{\rm (H)}$
have similar structures to $F_i^{\rm (K)}$, $F_i^{\rm (KW)}$, and $F_i^{\rm (H)}$, respectively.
%On the other hand, 
By contrast,
in view of \eqref{TwoTen}, we find that 
${\mathcal G}_{k,a}^{\rm (W)}$ may not be written in terms of the null forms,
differently from $F_i^{\rm (W)}$ with $N_1+1\le i \le N$.
Hence we divide $w^{\rm (ii)}$ into two parts $w^{\rm (iii)}$ and $w^{\rm (iv)}$:
For $1\le l \le N_4$, let $w_l^{\rm (iii)}$ and $w_l^{\rm (iv)}$ be the solutions
to
$$
\begin{cases}
\displaystyle \dal w_l^{\rm (iii)}=\sum_{a=0}^3 \pa_a
  \left\{{\mathcal G}_{N_3+l,a}(u, \pa u)-{\mathcal G}_{N_3+l,a}^{\rm (W)}(w, \pa w)\right\},\\
w_l^{\rm (iii)}(0,x)=w_l^{\rm (ii)}(0,x),\ 
\pa_t w_l^{\rm (iii)} (0,x)=\pa_t w_l^{\rm (ii)}(0,x),
\end{cases}
$$
and
$$
\begin{cases}
\displaystyle \dal w_l^{\rm (iv)}=F_{N_1+N_3+l}^{\rm(W)}(w, \pa w, \pa_x\pa w),\\
w_l^{\rm (iv)}(0,x)=\pa_t w_l^{\rm (iv)} (0,x)=0,
\end{cases}
$$
respectively.
Since we have
$$
F_{N_1+N_3+l}^{\rm (W)}(w,\pa w, \pa_x\pa w)=\sum_{a=0}^3 \pa_a\left({\mathcal G}_{N_3+l, a}^{\rm (W)}(w, \pa w) \right),\quad 1 \le l\le N_4,
$$
we get $w_l^{\rm (ii)}(=w_{N_3+l})=w_l^{\rm (iii)}+w_l^{\rm (iv)}$.
We put $w^{\rm (iii)}=(w^{\rm (iii)}_l)$ and $w^{\rm (\rm iv)}=(w^{\rm (iv)}_l)$ with $1\le l \le N_4$.

% Note that we have ${\mathcal G}_{k,a}-{\mathcal G}_{k,a}^{\rm (W)}
% ={\mathcal G}^{\rm (K)}_{k,a}+{\mathcal G}^{\rm (KW)}_{k,a}+{\mathcal G}^{\rm (H)}_{k,a}$.
It is easy to see that ${\mathcal G}_{k,a}^{\rm (K)}$, ${\mathcal G}_{k,a}^{\rm (KW)}$, and ${\mathcal G}_{k,a}^{\rm (H)}$ enjoy the estimates corresponding to \eqref{Right12}, \eqref{Right13} and \eqref{Right14}, respectively. 
% we get
% \begin{align}
% \label{Right12'}
% \bigl|{\mathcal G}_{k,a}^{\rm (K)}(v, \pa v) \bigr|_{s}\le & CA\ve\jb{t+|x|}^{-3/2}|v|_{s+1},\\
% \nonumber
% \bigl|{\mathcal G}_{k,a}^{\rm (KW)}(u, \pa u)\bigr|_s
% \le & CA\ve \jb{t+|x|}^{-1}W_-(t, |x|)^{-1+p} |v|_{s+1}\\
% &{}+CA\ve\jb{t+|x|}^{-3/2}\left(|w^{\rm (ii)}|_s+|\pa w|_{s}\right),
% \label{Right13'}\\
% \label{Right14'}
% \bigl|{\mathcal G}_{k,a}^{\rm (H)}(u,\pa u)\bigr|_s \le & CA^2\ve^2 \jb{t+|x|}^{-2+2\delta}\jb{t-|x|}^{-2\delta}\left(|v|_{s+1}+|w|_s+|\pa w|_{s}\right)
% \end{align}
% at $(t,x)\in [0, T)\times \R^3$, provided that $s\le 2\sigma$.
Since we have 
$$
{\mathcal G}_{k,a}-{\mathcal G}_{k,a}^{\rm (W)}
={\mathcal G}^{\rm (K)}_{k,a}+{\mathcal G}^{\rm (KW)}_{k,a}+{\mathcal G}^{\rm (H)}_{k,a},
$$
similarly to \eqref{KY10} we get
\begin{equation}
\jb{x}\bigl|{\mathcal G}_{k,a}-{\mathcal G}_{k,a}^{\rm (W)}\bigr|_{2\sigma-7}
\le C_*A\ve^2\jb{t+|x|}^{-1+3\lambda-\mu}W_-(t,|x|)^{-1+\mu},
\end{equation}
and Lemma~\ref{DivDecay} yields
\begin{equation}
\label{Shera}
\sup_{(t,x)\in [0,T)\times \R^3} \jb{t+|x|}^{-3\lambda}\jb{x}\jb{t-|x|}
|w^{\rm (iii)}(t,x)|_{2\sigma-8}\le C_*\ve.
\end{equation}
%On the other hand, in 
In view of \eqref{Right21}, using \eqref{KubYok} with $(\rho, \kappa)=
(0, (3/2)-2\lambda)$, we get
\begin{equation}
\label{Conie}
\sup_{(t,x)\in [0,T)\times \R^3} \jb{t+|x|}\jb{t-|x|}^{(1/2)-2\lambda}
|w^{\rm (iv)}(t,x)|_{2\sigma-8}\le C_*\ve,
\end{equation}
since we may assume $\mu+\delta-\lambda<0$ for small $\mu>0$.
Summing up the estimates, we obtain
\begin{equation}
\label{Decay12'}
\sup_{(t,x)\in [0,T)\times \R^3} \jb{t+|x|}^{1-3\lambda}W_-(t,|x|)^{(1/2)+\lambda}
|w^{\rm (ii)}(t,x)|_{2\sigma-8}\le C_*\ve,
\end{equation}
provided that $\ve$ is sufficiently small.
%%%%%%%%%%%%%%%%%%%%%%%%%%%%%%

{\it Step 4: Decay Estimates, Part 3.}
Let $1\le i\le N_1$ in this step.
Motivated by the technique in \cite{Bac88}, \cite{Kos92} and \cite{YTsu03},
we introduce
\begin{equation}
\label{YTran} \widetilde{v}_i=v_i-m_i^{-2}F_i^{\rm (W)}(w, \pa w, \pa_x\pa w)
\end{equation}
in order to treat $F_i^{\rm (W)}$ for which the null condition is not assumed. 
Then we get
\begin{align}
\nonumber
\left(\dal+m_i^2\right)\widetilde{v}_i=& \left\{\left(\dal+m_i^2\right)v_i
{}-F_i^{\rm (W)}(w, \pa w, \pa_x\pa w)\right\}
\\
&
{}-m_i^{-2}\dal \left\{F_i^{\rm (W)}(w, \pa w, \pa_x\pa w)\right\}.
\label{YTran01}
\end{align}
From the condition (b--i), we can write
\begin{equation}
\label{ExpFW}
F_i^{\rm (W)}(w,\pa w, \pa_x\pa w)=\sum_{|\alpha|,|\beta|\le 2} \sum_{1\le j, k\le N_2}
P_i^{jk\alpha\beta} (\pa^\alpha w_j)(\pa^\beta w_k)
\end{equation}
with appropriate constants $P_i^{jk\alpha\beta}$, where $P_i^{jk\alpha\beta}$ vanishes
either when $1\le j\le N_3$ and $|\alpha|=0$, or when $1\le k\le N_3$ and $|\beta|=0$.
Since we have 
$$
\dal(\varphi\psi)=2Q_0(\varphi, \psi)+(\dal\varphi)\psi+\varphi(\dal\psi)
$$
for any smooth functions $\varphi$ and $\psi$, from \eqref{ExpFW} we get
\begin{equation}
\label{YTran02}
\dal F_i^{\rm (W)}=\widetilde{F}_i^{\rm (W)}+\widetilde{F}_i^{\rm (H)},
\end{equation}
where
\begin{align}
\label{DefFTilde}
 \widetilde{F}_i^{\rm (W)}=& 2\sum_{|\alpha|, |\beta|\le 2}\sum_{1\le j,k\le N_2}
 P_i^{jk\alpha\beta} Q_0(\pa^\alpha w_j, \pa^\beta w_k),\\
\label{DefHTilde}
 \widetilde{F}_i^{\rm (H)}=& \sum_{|\alpha|, |\beta|\le 2}\sum_{1\le j,k\le N_2}
 P_i^{jk\alpha\beta}\left\{
                             (\pa^\alpha F_{N_1+j})(\pa^\beta w_k)+(\pa^\alpha w_j)(\pa^\beta F_{N_1+k})
                          \right\}.
\end{align}
Note that each $\widetilde{F}_i^{\rm (W)}$ 
is written in terms of the null forms, and we can expect
extra decay for $\widetilde{F}_i^{\rm (W)}$. 
%On the other hand, each 
Note also that each
$\widetilde{F}_i^{\rm (H)}$ is a function
of cubic order with respect to $\pa^\alpha u_j$ with $|\alpha|\le 4$ and $1\le j\le N$.
From \eqref{YTran01} and \eqref{YTran02}, we obtain
\begin{align}
\label{YTran03}
\left(\dal+m_i^2\right) \widetilde{v}_i=F_i^{\rm (K)}+F_i^{\rm (KW)}-m_i^{-2}\widetilde{F}_i^{\rm (W)}
{}+\bigl(F_i^{\rm (H)}-m_i^{-2}\widetilde{F}_i^{\rm (H)}\bigr).
\end{align}

By \eqref{Right12} and \eqref{Ene01}, we have
\begin{align}
\nonumber
\left\|
 \jb{t+|\cdot|}
  \bigl| F_i^{\rm (K)} \bigr|_{2\sigma-10}
\right\|_{L^2(\R^3)}\le & CA\ve (1+t)^{-1/2}\|v\|_{2\sigma-8}
\\
\le & C_* A\ve^2 (1+t)^{\lambda-(1/2)}.
\label{Right32}
\end{align}

Similarly to \eqref{Right21} and \eqref{Right23},
but using \eqref{Decay11}, \eqref{Decay12} and \eqref{Decay12'} instead of \eqref{Decay02}, \eqref{Decay03} and
\eqref{Decay03'}, we obtain
\begin{align}
\label{Right31}
\jb{t+|x|}\bigl| \widetilde{F}_i^{\rm (W)} \bigr|_{2\sigma-10}\le & C_*A\ve^2 \jb{x}^{-1}\jb{t+|x|}^{-1+\delta+3\lambda}{W}_-(t,|x|)^{-1-\delta},\\
\label{Right33}
\jb{t+|x|}\bigl| F_i^{\rm (KW)} \bigr|_{2\sigma-10}\le & C_*A\ve^2 \jb{t+|x|}^{-(3/2)+3\lambda}{W}_-(t,|x|)^{-(1/2)-\lambda}
\end{align}
at $(t,x)\in [0,\infty)\times \R^3$.
%(for the later usage, we note that \eqref{Right33} is also true for $1\le i\le N$).
Since we may assume $\delta<1/2$, it follows from \eqref{Right31} and \eqref{Right33}
that
\begin{align}
\nonumber
& \left\|
       \jb{t+|\cdot|}
        \left(
          \bigl| \widetilde{F}_i^{\rm (W)}\bigr|_{2\sigma-10}
           {}+\bigl|F_i^{\rm (KW)}\bigr|_{2\sigma-10}
        \right)
    \right\|_{L^2(\R^3)}\\
& \qquad 
\le C_*A\ve^2 \left\|\jb{\, \cdot\, }^{-1}\jb{t+|\cdot|}^{-(1/2)+3\lambda} {W}_-(t,|\cdot|)^{-(1/2)-\lambda}\right\|_{L^2(\R^3)}
\nonumber\\
& \qquad \le C_*A\ve^2(1+t)^{-(1/2)+3\lambda}.
\label{Right41}
\end{align}

Going similar lines to \eqref{Right14}, and then using \eqref{Ene01} and \eqref{Decay11},
we obtain
\begin{align}
\nonumber
& \left\|
     \jb{t+|\cdot|} \bigl|F_i^{\rm (H)}+m_i^{-2} \widetilde{F}_i^{\rm (H)}\bigr|_{2\sigma-10}
   \right\|_{L^2(\R^3)}\\
& \qquad \le CA^2\ve^3 (1+t)^{-1+2\delta}
\left(\|v\|_{2\sigma-10}+\|w^{\rm (ii)}\|_{2\sigma-10}+\|\pa u\|_{2\sigma-7}\right) \nonumber\\
& \quad\qquad {}+C_*A^2\ve^3\left\|\jb{t+|\cdot|}^{-2+3\delta+3\lambda}\jb{t-|\cdot|}^{-3\delta}
\right\|_{L^2(\R^3)}
\nonumber\\
& \qquad \le C_*A^2\ve^3\left((1+t)^{-1+\lambda+2\delta}+(1+t)^{-(1/2)+3\lambda}\right).
\end{align}
To sum up, we have proved
\begin{equation}
\left\|\jb{t+|\cdot|}
  \left|\left(\dal+m_i^2\right) \widetilde{v}_i\right|_{2\sigma-10}
\right\|_{L^2} %\le & C_*A\ve^2 (1+t)^{-(1/2)+3\lambda} 
\le  
C_*A\ve^2(1+t)^{-1/4}, 
\end{equation}
because $\lambda\le 1/12$ and $\delta\ll 1$.
Now Lemma \ref{DecayKleinGordon} implies
\begin{align}
\nonumber
& \sup_{(t,x)\in [0, T)\times \R^3}\jb{t+|x|}^{3/2}|\widetilde{v}(t,x)|_{2\sigma-14}
\\
& \quad\qquad 
\le C_*\Bigl(\ve+A\ve^2\sum_{j=0}^\infty \sup_{\tau\in (0,t)}\chi_j(\tau)(1+\tau)^{-1/4}\Bigr)
\le C_*\ve,
\label{Decay41-o}
\end{align}
because we have $\sup_{\tau\in (0,t)}\chi_j(\tau)(1+\tau)^{-1/4}\le 2^{-(j-1)/4}$ for $j\ge 1$.
From \eqref{Decay12}, \eqref{Decay12'} and \eqref{ExpFW}, we get
\begin{align*}
\bigl|F_i^{\rm (W)}\bigr|_{2\sigma-14} \le & C\left(|\pa w|_{\sigma+1}+|w^{\rm (ii)}|_{\sigma+2}\right)
\left(|\pa w|_{2\sigma-13}+|w^{\rm (ii)}|_{2\sigma-14}\right)\\ 
\le & C_*A\ve^2\jb{t+|x|}^{3\lambda-2}, \nonumber
\end{align*}
which, together with \eqref{YTran} and \eqref{Decay41-o}, yields
\begin{equation}
\label{Decay41}
\sup_{(t,x)\in [0, T)\times \R^3}\jb{t+|x|}^{3/2}|v(t,x)|_{2\sigma-14}\le C_*\ve 
\end{equation}
for small $\ve$.
%%%%%%%%%%%%%%%%%%%%%%%

{\it Step 5: Decay Estimates, the Final Part.} Let $N_1+1\le i\le N$.
%%%%%%%%%%%%%%%%%%%%%%%%%%%%%%%
Similarly to \eqref{Right31}, we get
\begin{equation}
\label{AA1}
\jb{x}\bigl|F_i^{\rm (W)}\bigr|_{2\sigma-9}\le C_*A\ve^2 \jb{t+|x|}^{-2+\delta+3\lambda}
W_-(t,|x|)^{-1-\delta}.
\end{equation}
Therefore, using \eqref{KubYok} with $(\rho, \kappa)=(0, 2-4\lambda)$, we get
\begin{equation}
\label{AA0}
\sup_{(t,x)\in [0, T)\times \R^3} \jb{t+|x|}\jb{t-|x|}^{1-4\lambda}
 \bigl|w^{\rm (iv)}(t,x)\bigr|_{2\sigma-9}\le C_*\ve,
\end{equation}
which, together with \eqref{Shera}, yields
\begin{equation}
\label{Celenaria}
\sup_{(t,x)\in [0, T)\times \R^3}\jb{t+|x|}^{1-3\lambda}W_-(t,|x|)^{1-\lambda}
\bigl|w^{\rm (ii)}(t,x)\bigr|_{2\sigma-9}\le C_*\ve.
\end{equation}
%%%%%%%%%%%%%%%%%%%%%%%%%%%%%%%
By \eqref{Right12} and \eqref{Decay41}, we obtain
\begin{equation}
\label{Right61}
\jb{x}\bigl|F_i^{\rm (K)}\bigr|_{2\sigma-16}\le C_*A\ve^2 \jb{t+|x|}^{-2}\le C_*A\ve^2\jb{t+|x|}^{-1-\mu}W_-(t, |x|)^{\mu-1}
\end{equation}
at $(t,x)\in [0, T)\times \R^3$ for small $\mu>0$. 
From \eqref{Right13}, \eqref{Decay12}, \eqref{Decay41}, and \eqref{Celenaria} we get
\begin{align}
\nonumber
\jb{x}\bigl|F_i^{\rm (KW)}\bigr|_{2\sigma-16}\le & 
C_*A\ve^2 \jb{t+|x|}^{3\lambda-(3/2)}W_-(t,|x|)^{-1+p}
\\
\le & C_*A\ve^2 \jb{t+|x|}^{-1-\mu}W_-(t, |x|)^{\mu-1}
\label{AA2}
\end{align}
for small $\mu>0$.
%On the other hand, 
\eqref{Right14}, \eqref{Decay12}, \eqref{Decay41},
and \eqref{Celenaria} yield
\begin{align}
\nonumber
\jb{x}\bigl|F_i^{\rm (H)}\bigr|_{2\sigma-16}\le & 
C_*A^2\ve^3\jb{t+|x|}^{-(5/2)+2\delta}\jb{t-|x|}^{-2\delta}\\
&{}+C_*A^2\ve^3\jb{t+|x|}^{-2+2\delta+3\lambda}W_-(t,|x|)^{\lambda-2\delta-1}
\nonumber\\
& {}+CA^2\ve^2M^{\rm (i)}_{0,2\sigma-16} \jb{t+|x|}^{-2+3\delta}\jb{t-|x|}^{-3\delta}\nonumber\\
\le & (C_*A^2\ve^3+C A^2\ve^2M^{\rm (i)}_{0, 2\sigma-16})\jb{t+|x|}^{-1-\mu}W_-(t, |x|)^{\mu-1}
\label{Right62}
\end{align}
at $(t,x)\in [0,T)\times \R^3$.
Since we may assume $\delta+3\lambda\le 1/4$, say,
it follows from \eqref{AA1}, \eqref{Right61}, \eqref{AA2}, and \eqref{Right62}
that
\begin{equation}
\jb{x}|F_i|_{2\sigma-16}\le (C_* A\ve^2+CA^2\ve^2M^{\rm (i)}_{0, 2\sigma-16})\jb{t+|x|}^{-1-\mu}{W}_-(t,|x|)^{\mu-1}.
\label{lalala2}
\end{equation}
Now Lemma \ref{Asakura} and \eqref{KubYok} with $(\rho, \kappa)=(0, 1)$ imply
$$
M^{\rm (i)}_{0,2\sigma-16}\le C_*(\ve+A\ve^2)+CA^2\ve^2 M^{\rm (i)}_{0, 2\sigma-16},
$$
which leads to
\begin{equation}
\label{Decay62}
\sup_{(t,x)\in[0,T)\times \R^3} \jb{t+|x|}{\mathcal W}_0(t,|x|)|w^{\rm (i)}(t,x)|_{2\sigma-16}
=M^{\rm (i)}_{0, 2\sigma-16} \le C_*\ve,
\end{equation}
provided that $\ve$ is small enough.
By \eqref{lalala2} and \eqref{Decay62}, using Lemma \ref{Asakura} and
\eqref{KubYokD} with $(\rho, \kappa)=(0,1)$, we obtain
\begin{equation}
\label{Decay63}
\sup_{(t,x)\in[0,T)\times \R^3} \jb{x}\jb{t-|x|}|\pa w(t,x)|_{2\sigma-17}
\le C_*\ve, 
\end{equation}
provided that $\ve$ is small enough.

%In view of \eqref{Right12'}, \eqref{Right13'}, and \eqref{Right14'},
Going similar lines to \eqref{Right61}--\eqref{lalala2}, and
using \eqref{Decay62}, we obtain
\begin{equation}
\jb{x} \bigl|{\mathcal G}_{k,a}-{\mathcal G}_{k,a}^{\rm (W)}\bigr|_{2\sigma-16}
\le C_*A\ve^2 \jb{t+|x|}^{-1-\mu} W_-(t,|x|)^{\mu-1},
\end{equation}
and Lemma~\ref{DivDecay} yields
\begin{equation}
\label{AA3}
\sup_{(t,x)\in[0,T)\times \R^3} \jb{x}\jb{t-|x|}|w^{\rm (iii)}(t,x)|_{2\sigma-17}
\le C_*\ve.
\end{equation}
Now, \eqref{AA0} and \eqref{AA3} imply
\begin{equation}
\label{Decay63'}
\sup_{(t,x)\in[0,T)\times \R^3} \jb{t+|x|}W_-(t,|x|)^{1-p}|w^{\rm (ii)}(t,x)|_{2\sigma-17}
\le C_*\ve, 
\end{equation}
provided that $\ve$ is small enough, since we have assumed $4\lambda< p$ at the beginning.

%%%%%%%%%%%%%%%%%%%%

{\it Step 6: Conclusion.} Finally, \eqref{Decay41}, \eqref{Decay62}, \eqref{Decay63}, and \eqref{Decay63'} yield
\begin{equation}
\sup_{0\le t<T} \|d_{\sigma, p}(t,\cdot)\|_{L^\infty(\R^3)} \le C_0 \ve
 \label{LastIneq}
\end{equation}
for $\ve\le \ve_0(A)$, where $\ve_0(A)$
is a positive constant depending on $A$, and $C_0$ is some positive constant
which depends on $M_0$ and $F$, 
but is independent of $A$, $\ve$ and $T$.
We put $A_0=2C_0$.
Now \eqref{LastIneq} implies \eqref{boot-finale} for $A\ge A_0$ and $\ve\le \ve_0(A)$.
This completes the proof of Proposition \ref{bootstrap}.
%%%%%%%%%%%%%%%%%%%%%%%%%%%%%%%%%%%
%%%%%%%%%%%%%%%%%%%%%%%%%%%%%%%%%%%
\section*{Acknowledgments}
The author is grateful to Professors Masahito Ohta and Hideaki Sunagawa
for giving the author valuable informations on the previous works in the related field.
The author would also like to thank Professor Hideo Kubo for fruitful conversations.
This research is partially supported by Grant-in-Aid for Scientific Research (C) (No.~20540211), Japan Society for the Promotion of Science.
%%%%%%%%%%%%%%%%%%%%%%%%
%%%%%%%%%%%%%%%%%%%%%%%%%%%%%%%%%%%%%%%%%%%%%%%%
%%%%%%%%%%%%%%%% References %%%%%%%%%%%%%%%%%%%%%
%%%%%%%%%%%%%%%%%%%%%%%%%%%%%%%%%%%%%%%%%%%%%%%%%
{\small

\begin{flushleft}
{\sc Department of Mathematics, Wakayama University\\
930 Sakaedani, Wakayama 640-8510, Japan} \\
e-mail: katayama@center.wakayama-u.ac.jp
\end{flushleft}
}

\begin{thebibliography}{99}
%%%%%%%%%%%%%%%%%%%%% AAA
%%%%%%%%%%%%%%%%%%%%%%%%%%%%%%%%%%
\bibitem{Asa86} F.~Asakura, 
{\it Existence of a global solution to a semi-linear wave equation
with slowly decreasing initial data in three space dimensions},
Comm.~Partial Differential Equations
{\bf 11}, 1459--1487 (1986).

\bibitem{Bac88} A.~Bachelot, 
 {\it Probl\`eme de Cauchy global pour des syst\`emes de Dirac-Klein-Gordon},
Ann.~Inst.~Henri Poincar\'e {\bf 48}, 387--422 (1988).

\bibitem{Chr86} D.~Christodoulou, 
{\it Global solutions of nonlinear hyperbolic equations for small initial data},
{Comm.~Pure Appl.~Math.} {\bf 39}, {267--282} (1986).


%%%%%%%%%%%%%%%% DDD

%%%%%%%%%%%%%%%% GGG
\bibitem{Geo90} V.~Georgiev, 
 {\it Global solution of the system of  wave and Klein-Gordon equations},
 Math.~Z. {\bf 203}, 683--698 (1990).
\bibitem{Geo92} V.~Georgiev, 
 {\it Decay estimates for the Klein-Gordon equation},
 Comm.~Partial Differential Equations {\bf 17}, 1111--1139 (1992).
%%%%%%%%%%%%%%%% HHH

\bibitem{Hay-Nau-Bag08}
N.~Hayashi, P.~I.~Naumkin, and Ratno Bagus Edy Wibowo,
 {\it Nonlinear scattering for a system of nonlinear Klein-Gordon equations},
 J.~Math.~Phys. {\bf 49}, 103501 (2008).
\bibitem{Hoe88} L.~H\"ormander, 
                {$ L^{\,1}, L^{\,\infty} $ estimates for the wave
                     operator},
                in: Analyse Math\' ematique et Applications,
                   Contributions en l'Honneur de J. L. Lions,
                pp.~211--234, Gauthier-Villars, Paris (1988). 
                
\bibitem{Hoe97} L.~H\" ormander,   
                {Lectures on Nonlinear Hyperbolic Differential Equations},
                 Springer-Verlag, Berlin (1997). 
%%
\bibitem{Hos-Kub00}
   A.~Hoshiga, and H.~Kubo,
   {\it Global small amplitude solutions
               of nonlinear hyperbolic systems with a critical
                   exponent under the null condition}, 
   SIAM J.~Math.~Anal. {\bf 31}, 486--513 (2000).

%%%%%%%%%%%%%%%%%%%%%  JJJ
\bibitem{Joh79} F.~John,  
    {\it Blow-up of solutions of nonlinear wave equations 
                      in three space dimensions},
   Manuscripta Math. {\bf 28}, 235--268 (1979).
\bibitem{Joh81} F.~John, 
{\it Blow-up of solutions for quasi-linear wave equations in three
space dimensions}, Comm.~Pure Appl.~Math. {\bf 34}, 29--51 (1981).

%%%%%%%%%%%%%%%%%%%%%  KKK

\bibitem{Kat04:02} S.~Katayama, 
    {\it Global and almost-global existence for
                systems of nonlinear wave equations with different
                propagation speeds}, 
  Diff.~Integral Eqs. {\bf 17}, 1043--1078 (2004).

\bibitem{Kat-Kub08:01} S.~Katayama, and H.~Kubo, 
{\it Decay estimates of a tangential derivative to 
the light cone for the wave equation and their application},
SIAM J.~Math.~Anal. {\bf 39}, 1851--1862 (2008).
\bibitem{Kat-Kub08:03} S.~Katayama, and H.~Kubo, 
{\it An alternative proof of global existence for nonlinear wave 
equations in an exterior domain},
J.~Math.~Soc.~Japan {\bf 60}, 1135--1170 (2008).
\bibitem{Kat-Yok06} S.~Katayama, and K.~Yokoyama, 
{\it Global small amplitude solutions to systems of nonlinear wave
equations with multiple speeds}, 
Osaka J.~Math. {\bf 43}, 283--326 (2006).

\bibitem{Kla85:02} S.~Klainerman, 
{\it Global existence of small amplitude solutions to nonlinear Klein-Gordon equations
in four space-time dimensions},
 Comm.~Pure Appl.~Math. {\bf 38}, 631--641 (1985).

\bibitem{Kla86} S.~Klainerman,  
{\it The null condition and global existence to nonlinear wave equations},
in: Nonlinear Systems of Partial Differential Equations in Applied Mathematics, Part 1, Lectures in Applied Math. {\bf 23}, pp.~{293--326}, AMS, Providence, RI (1986).

\bibitem{Kla87} S.~Klainerman,  
  {\it Remarks on the global Sobolev inequalities
                      in the Minkowski space {${\mathbf R}^{n+1}$}},
  Comm.~Pure Appl. Math. {\bf 40}, 111--117 (1987).

\bibitem{Kos92} R.~Kosecki, 
          {\it The unit condition and global existence for a class of nonlinear Klein-Gordon equations},
           J.~Differential Equations {\bf 100}, 257--268 (1992).
\bibitem{Kov87} 
   M.~Kovalyov, 
           {\it Long-time behavior of solutions of a system of
                                          nonlinear wave equations},
             Comm.~Partial Differential Equations
             {\bf 12}, 471--501 (1987).
%%%%%%%%%%%%%%%%%%%%%%%%%%%%%%%%%%%%%

\bibitem{Kov-Tsu} M.~Kovalyov, and K.~Tsutaya, 
{\it Erratum: ``Long-time behavior of solutions of a system of
                  nonlinear wave equations''},
Comm.~Partial Differential Equations {\bf 18}, 1971--1976 (1993).

\bibitem{Kub-Yok01} K.~Kubota, and K.~Yokoyama,  
     {\it Global existence of classical solutions to
                       systems of nonlinear wave equations with different
                       speeds of propagation},
  Japanese J.~Math. {\bf 27}, 113--202 (2001).
  
%%%%%%%%%%%%%%%%%%%%%%%%%% LLL


%%%%%%%%%%%%%%%%%%%%%%%%%% MMM

%%%%%%%%%%%%%%%%%%%%%%%%%%%% N %%%%%%%%%%%%%%%%%%%%%%%%
%%%%%%%%%%%%%%%%%%%%%%%%%%%% O %%%%%%%%%%%%%%%%%%%%%%%%
\bibitem{OzaTsuTsu95} T.~Ozawa, K.~Tsutaya, and Y.~Tsutsumi, 
 {\it Normal form and global solutions for the Klein-Gordon-Zakharov equations},
 Ann.~Inst.~H.~Poincar{\'e} Anal.~Non Lin{\'e}aire {\bf 12}, 459--503 (1995).
\bibitem{OzaTsuTsu99} T.~Ozawa, K.~Tsutaya, and Y.~Tsutsumi, 
 {\it Well-posedness in energy space for the Cauchy problem of the Klein-Gordon-Zakharov
 equations with different propagation speeds in three space dimensions},
 Math.~Ann. {\bf 313}, 127--140 (1999).
%%%%%%%%%%%%%%%%%%%%%%%%%% SSS

\bibitem{Sha85} J.~Shatah, 
 {\it Normal forms and quadratic nonlinear Klein-Gordon Equations},
 Comm.~Pure Appl. Math. {\bf 38}, 685--696 (1985).


%%%%%%%%%%%%%%%%%%%%%%%%%%%%%%%%%%%%%%%%
\bibitem{Sid89} T.~C.~Sideris, 
                 {\it Decay estimates for the three space dimensional inhomogeneous 
                  Klein-Gordon equation and applications},
                 Comm.~Partial Differential Equations {\bf 14}, 1421--1455 (1989).

\bibitem{Str68} W.~A.~Strauss, 
                {\it Decay and asymptotics for $\dal u=F(u)$},
                 J.~Funct.~Anal. {\bf 2}, 409--457 (1968).
%%%%%%%%%%%%%%%%%%%%%%%%%% TTT
\bibitem{KTsu96} K.~Tsutaya, 
                {\it Global existence of small amplitude solutions for the 
                Klein-Gordon-Zakharov equations},
                Nonlinear Anal.~{\bf 27}, 1373--1380 (1996).
\bibitem{YTsu03} Y.~Tsutsumi, 
                {\it Global solutions for the Dirac-Proca equations with
                small initial data in $3+1$ space time dimensions},
                J.~Math.~Anal.~Appl. {\bf 278}, 485--499 (2003).
%%%%%%%%%%%%%%%%%%%%%%%%%% YYY
\bibitem{Yok00} K.~Yokoyama, 
   {\it Global existence of classical solutions to 
              systems of wave equations with critical nonlinearity
              in three space dimensions},
  J.~Math.~Soc.~Japan {\bf 52}, 609--632 (2000).
\end{thebibliography}
\end{document}